\input ./style/arxiv-general.cfg
\documentclass[aop,MSNbibl,seceqn,dvips]{arximspdf}
\makeatletter
   \@ifpackageloaded{graphicx}{}{\usepackage{graphicx}}
\makeatother
\usepackage{mathbh}

\doi{10.1214/15-AOP1042}
\volume{44}
\issue{4}
\pubyear{2016}
\firstpage{3076}
\lastpage{3110}
\docsubty{FLA}

\makeatletter
\def\sfrac#1#2{#1/#2}
\def\vfrac#1#2{(#1)/#2}
\def\afrac#1#2{#1/(#2)}

\def\sklfrac#1#2{(#1/#2)}
\def\sklvfrac#1#2{((#1)/#2)}
\def\sklafrac#1#2{(#1/(#2))}

\newcommand{\rrvert}{\vert}

\newcommand{\llvert}{\vert}
\renewcommand{\mid}{|}
\newcommand{\mathds}{\mathbh}
\newtheorem{theorem}{Theorem}[section]
\newtheorem{lemma}[theorem]{Lemma}
\newtheorem{proposition}[theorem]{Proposition}
\newtheorem{corollary}[theorem]{Corollary}
\newproclaim{rem}[theorem]{Remark}
\newproclaim{notation}[theorem]{Notation}
\newproclaim{definition}[theorem]{Definition}
\newcommand{\R}{\mathbb{R}}
\newcommand{\Z}{\mathbb{Z}}
\newcommand{\Q}{\mathbb{Q}}
\newcommand{\E}{\mathds{E}}
\newcommand{\Pb}{\mathds{P}}
\newcommand{\ind}{\mathbh{1}}
\def\eps{\varepsilon}
\def\S{\mathbb{S}}
\def\T{\mathbb{T}}
\def\LB{\mathcal{B}} 
\def\LZ{\mathcal{Z}}
\def\u{\mathbf{u}}
\makeatother

\begin{document}
\begin{frontmatter}

\title{Liouville Brownian motion}
\runtitle{Liouville Brownian motion}

\begin{aug}
\author[A]{\fnms{Christophe}~\snm{Garban}\corref{}\thanksref{T1}\ead[label=e1]{garban@math.univ-lyon1.fr}},
\author[B]{\fnms{R\'emi}~\snm{Rhodes}\thanksref{T2}\ead[label=e2]{remi.rhodes@u-pem.fr}}
\and
\author[C]{\fnms{Vincent}~\snm{Vargas}\ead[label=e3]{Vincent.Vargas@ens.fr}\thanksref{T2}}
\runauthor{C. Garban, R. Rhodes and V. Vargas}
\affiliation{Universit\'e Lyon 1,
Universit{\'e} Paris-Est, Marne la Vall\'ee
and CNRS, ENS Paris}
\address[A]{C. Garban\\
Institut Camille Jordan\\
Universit\'e Lyon 1\\
43 bd du 11 novembre 1918\\
69622 Villeurbanne cedex \\
France \\
\printead{e1}}
\address[B]{R. Rhodes\\
LAMA\\
Universit{\'e} Paris-Est Marne la Vall\'ee\\
Champs sur Marne\\
France\\
\printead{e2}}
\address[C]{V. Vargas\\
DMA\\
ENS Paris\\
45 rue d'Ulm\\
75005 Paris\\
France \\
\printead{e3}}
\end{aug}
\thankstext{T1}{Supported in part by the ANR Grant MAC2 10-BLAN-0123.}
\thankstext{T2}{Supported in part by Grant ANR-11-JCJC CHAMU.}

%
\received{\smonth{5} \syear{2014}}
%
\revised{\smonth{6} \syear{2015}}

%
\begin{abstract}
We construct a stochastic process, called the \textit{Liouville Brownian
motion},
which is the Brownian motion associated to the metric $e^{\gamma
X(z)}\,dz^2$,
$\gamma<\gamma_c=2$ and $X$ is a Gaussian Free Field. Such a process
is
conjectured to be related to the scaling limit of random walks on large
planar maps eventually weighted by a model of statistical physics which
are embedded in the Euclidean plane or in the sphere in a conformal manner.
The construction amounts to changing the speed of a standard
two-dimensional
Brownian motion $B_t$ depending on the local behavior of the Liouville
measure
``$M_\gamma(dz) = e^{\gamma X(z)} \,dz$''. We prove that the associated
Markov
process is a Feller diffusion for all $\gamma<\gamma_c=2$ and that for
all
$\gamma<\gamma_c$, the Liouville measure $M_\gamma$ is invariant under
$P_\mathbf{t}$.
This Liouville Brownian motion enables us to introduce a whole set of
tools of
stochastic analysis in Liouville quantum gravity, which will be
hopefully useful in analyzing the geometry of Liouville quantum gravity.
\end{abstract}

%
\end{frontmatter}

\section{Introduction}

An important issue for applications in $2d$-Liouville quantum gravity
is to construct a random metric on a two-dimensional Riemann manifold
$D$, say a domain of $\R^2$ (or the sphere) equipped with the Euclidean
metric $dz^2$, which takes on the form
%
\begin{equation}
\label{i.metric} e^{\gamma X(z)}\,dz^2,
\end{equation}
where $X$ is a Gaussian Free Field (GFF) on the manifold $D$ and
$\gamma
\in[0,2)$ is a coupling constant that can be expressed in terms of the
central charge of the model coupled to gravity (see \cite{cfKPZ,cfDa}
for further details and also \cite{DistKa,GM,Nak} for insights in
Liouville quantum gravity). The simplicity of such an expression hides
many highly nontrivial mathematical difficulties. Indeed, the
correlation function of a GFF presents a short scale logarithmically
divergent behavior that makes relation~(\ref{i.metric}) nonrigorous.
One has to apply a cutoff procedure to smooth down the singularity of
the GFF and the method to do this at a metric level remains unclear.
However, many geometric quantities are related to this metric and for
some of them, the cutoff procedure may be applied properly. For
instance, the construction of the volume form, also called the
Liouville measure, was carried out by Kahane within the framework of
Gaussian multiplicative chaos \cite{cfKah} (see also \cite
{cfRoVa,cfRoVa1,cfDuSh} for more recent constructions based on
convolution techniques). This allows us to give a rigorous meaning to
the expression
%
\begin{equation}
\label{i.measure} M_\gamma(A)=\int_A
e^{\gamma X(z)-\sklfrac{\gamma^2}{2}\E[X(z)^2]} \,dz,
\end{equation}
where $dz$ stands for the volume form (Lebesgue measure) on $D$ [to be
exhaustive, one should integrate against $h(z) \,dz$ where $h$ is a
deterministic function involving the conformal radius at $z$ but this
term does not play an important role for our concerns]. This strategy
made possible an interpretation in terms of measures of the
Knizhnik--Polyakov--Zamolodchikov formula (KPZ for short, see \cite
{cfKPZ}) relating the fractal dimensions of sets as seen by the
Lebesgue measure or the Liouville measure. The KPZ formula is proved in
\cite{cfDuSh} when considering the fractal notion of expected box
counting dimension whereas the fractal notion of almost sure Hausdorff
dimension is considered in \cite{Rnew10,Rnew4} (see also \cite{Benj}).
The reader may consult \cite
{BKNSW,Benj,cfDa,David-KPZ,DistKa,Rnew7,Rnew12,bourbaki,GM,Nak,review}
for more references on this topic. Another important part of the theory
which we do not review here is that it is conjectured to be the scaling
limit of discrete quantum gravity: the reader may consult \cite
{cfDuSh} for more on this topic as well as physics references therein.

Another powerful tool in describing a Riemann geometry is the Brownian
motion. With it are attached several analytic objects serving to
describe the geometry: a semigroup, a Laplace--Beltrami operator, a
heat kernel, Dirichlet forms, etc. Therefore, a relevant way to have
further insights into Liouville quantum gravity geometry is to define
the \textit{Liouville Brownian motion} (LBM for short). This is the
purpose of this paper. It can be constructed on any background
$2d$-Riemann manifold equipped with a GFF $X$ and can be seen as the
Brownian motion associated to the metric $e^{\gamma X(x)} \,dx^2$ where
$\gamma\in[0,2[$ is a parameter and $dx^2$ stands for the metric on the
manifold.

In this paper and for pedagogical purposes, we will mostly describe the
situation when the underlying manifold is the whole plane $\R^2$, in
which case it is natural to consider a Massive Gaussian Free Field $X$
on $\R^2$ (MFF for short). We will also explain how to adapt our
framework to the cases of the sphere $\S^2$, the torus $\T^2$ or planar
bounded domains. More generally, it is also clear that our methodology
may apply to any $2$-dimensional Riemann manifold equipped with a
log-correlated Gaussian field and yields similar results. Let us also
mention that another work \cite{berest} appeared online simultaneously
to ours and is concerned with the LBM starting from one point: the
paper \cite{berest} proves that for fixed $x$, one can define almost
surely~the LBM starting from $x$ (the work \cite{berest} also initiates
a multi-fractal analysis of the LBM starting from one point). We will
show that almost surely we can define the LBM starting from all $x$,
hence obtaining the existence of a diffusion process associated to the
tensor~(\ref{i.metric}) and all the related stochastic analysis tools.

Finally, we point out that the notions of diffusion or heat kernel are
at the core of the physics literature about Liouville quantum gravity
(see \cite{amb,amb2,calg1,calg2,david-hd,David-KPZ,wata}, e.g.,
among a huge amount of other works). For instance, a heat kernel
derivation of the KPZ formula is obtained in \cite{David-KPZ}. The
fractal structure of quantum space--time is also investigated in \cite
{amb,amb2,david-hd,wata} via diffusions and heat kernel properties,
obtaining relations about the fractal dimensions of quantum space--time.

\section{Liouville Brownian motion on the plane}\label{sec.massLBM}
In this section, we construct Liouville Brownian motion on the whole
plane. Regarding the physics literature, the natural free field to
consider on the whole plane is the Massive Gaussian Free Field (MFF for
short). So we first remind the reader of the construction of the MFF
after introducing a few basic notation. Then we recall the construction
of Gaussian multiplicative chaos associated to the MFF and state a few
basic properties, which are used thereafter to construct the Liouville
Brownian motion.

\subsection{Basic notation and terminology}\label{index}

\subsubsection*{Basics}
In what follows, the Liouville Brownian motion that we are going to
construct will be denoted by $(\LB_\mathbf{t})_{\mathbf{t}\geq0}$.
We distinguish the
(quantum) time $\mathbf{t}$ along $\LB_\mathbf{t}$ and the
(classical) time $t$ along a
standard Brownian motion $B_t$.
The open ball centered at $x$ with radius $R$ is denoted $B(x,R)$ and
the closed ball~$\bar{B}(x,R)$. $\S^2$ and $\T^2$ stand for the
two-dimensional sphere or torus.

\subsubsection*{Functional spaces and analysis}
The space of continuous functions with compact support in a domain $D$
(resp., vanishing at infinity on $\R^2$, resp., bounded functions on
$D$, resp., continuous functions on $\R_+$ equipped with the sup-norm
topology over compact sets) is denoted by $C_c(D)$ [resp., $C_0(\R^2)$,
resp., $C_b(D)$, resp., $C(\R_+)$].

The standard Laplace--Beltrami operator on a manifold is denoted by
$\Delta$. We will say that a semigroup on $C_b(D)$ is Feller if the
semigroup maps $C_b(D)$ into itself. We will say that a Markov process
on $D$ is Feller if its semigroup is.

\subsubsection*{Positive continuous additive functionals and Revuz measures}
Let us consider a standard Brownian motion $(\Omega_B,(B_t)_{t \geq
0},(\mathcal{F}_t)_{t \geq0},(\Pb_x^B)_{x \in D})$ in $D$ (with
$D=\R
^2$, $\S^2$ or $\T^2$). It is reversible for the canonical volume form
$dx$ of $D$. We suppose that the space $\Omega_B$ is equipped with the
standard shifts $(\theta_t)_{t \geq0}$ on the trajectory. One may then
consider the classical notion of capacity $\operatorname{Cap}$ associated to the
Brownian motion (see \cite{fuku}). The set $K$ is said polar when
$\operatorname{Cap}(K)=0$.

A Revuz measure $\mu$ is a Radon measure on $D$ which does not charge
the polar sets. A positive continuous additive functional (PCAF)
$(A_t)_{t \geq0}$ is a $\mathcal{F}_t$-adapted continuous functional
with values in $[0,\infty]$ that satisfies for all $\omega\in\Lambda$:
\[
A_{t+s}(\omega)=A_s(\omega)+ A_t\bigl(
\theta_s(\omega)\bigr), \qquad s,t \geq0,
\]
where $\Lambda$ is a subset of $\Omega_B$ such that $\forall x \in D$,
$\Pb_x^B(\Lambda)=1$ and $\theta_t(\Lambda) \subset\Lambda$ for
all $t
\geq0$. In particular, we will always consider PCAFs defined for all
starting points $x \in D$ (they are sometimes called \textit{PCAF in the
strict sense} in the literature, especially in \cite{fuku}).

\subsubsection*{Massive Gaussian Free Field on the plane}
We consider a whole plane Massive Gaussian Free Field (MFF) (see \cite
{glimm,She07} for an overview of the construction of the MFF and
applications). Given a real number $m>0$, it is a centered Gaussian
random distribution (in the sense of Schwartz) with covariance function
given by the Green function $G_m$ of the operator $m^2-\triangle$,
that is,
\[
\bigl(m^2-\triangle\bigr)G_m(x,\cdot)=2\pi
\delta_x.
\]
Notice that $G_m$ is $\pi$ times the Green function of the Brownian
motion killed at rate $m^2/2$.
It is a standard fact that the massive Green function can be written as
an integral of the transition densities of the Brownian motion weighted
by the exponential of the mass:
%
\begin{equation}
\label{MFF1} \forall x,y \in\R^2,\qquad G_m(x,y)=\int
_0^{\infty}e^{-\sklfrac {m^2}{2}u-\afrac{\llvert  x-y\rrvert  ^2}{2u}}\frac{du}{2 u}.
\end{equation}
Clearly, it is a kernel of $\sigma$-positive type in the sense of
Kahane \cite{cfKah} since we integrate a continuous function of
positive type with respect to a positive measure. One can also check that
%
\begin{equation}
\label{logkernel} G_m(x,y)=\ln_+\frac{1}{\llvert  x-y\rrvert  }+g_m(x,y),
\end{equation}
for some continuous and bounded function $g_m$ and $\ln_+x=\max(\ln
x, 0)$.

It is furthermore a star-scale invariant kernel (see \cite
{Rnew1,sohier}): it can be rewritten as
%
\begin{equation}
\label{MFF3} G_m(x,y)=\int_{1}^{+\infty}
\frac{k_m(u(x-y))}{u} \,du,
\end{equation}
for some continuous covariance kernel $k_m(z)=\frac{1}{2}\int_0^\infty
e^{-\sklafrac{m^2}{2v}\llvert  z\rrvert  ^2-\sfrac{v}{2}} \,dv$. In particular, we will make
intensive use of the following relation, valid for $\epsilon\in\,]0,1]$:
%
\begin{equation}
\label{ssi} G_m(x,y)\leq G_m\biggl(
\frac{x}{\epsilon},\frac{y}{\epsilon}\biggr)+\ln\frac
{1}{\epsilon}.
\end{equation}

Now we consider an unbounded strictly increasing sequence $(c_n)_{n\geq
1}$ such that $c_1=1$. For each $n\geq1$, we consider a centered
Gaussian process $Y_n$ with covariance kernel given by
%
\begin{equation}
\label{kernelY} \E\bigl[Y_n(x)Y_n(y)\bigr]=\int
_{c_{n-1}}^{c_{n}}\frac{k_m(u(x-y))}{u} \,du.
\end{equation}
The reader may check that such a process is stationary and has smooth
sample paths (to check this point, apply the Kolmogorov criterion to
$Y_n$ as well as its derivatives in the standard manner).
The MGFF is the Gaussian distribution defined by
\[
X(x)=\sum_{n\geq1}Y_n(x),
\]
where the processes $(Y_n)_n$ are assumed to be independent.
We define the $n$-re\-gularized field by
%
\begin{equation}
\label{MGFFdec} X_n(x)=\sum_{k= 1}^nY_k(x).
\end{equation}
Actually, based on Kahane's theory of multiplicative chaos \cite
{cfKah}, the choice of the decomposition~(\ref{MGFFdec}) will not play
a part in the forthcoming results, except that it is important that the
covariance kernel of $X_n$ be smooth in order to associate to this
field a Riemann geometry.

\begin{notation}\label{product} In what follows, we will consider a
Brownian motion on $\R^2$ (or other two-dimensional manifolds)
$(\Omega
_B,\mathcal{F}_B,(B_t)_{t \geq0},(\mathcal{F}_t)_{t \geq0},(\Pb
^B_x)_{x \in\R^2})$. We will also consider a MFF $X$ (and all the
corresponding $(Y_n)_n$) defined on a probability space $(\Omega
_X,\mathcal{F}_X,\Pb^X)$. So we consider a measurable space $(\Omega
,\mathcal{F})=(\Omega_X\times\Omega_B,\mathcal{F_X\otimes\mathcal
{F}_B})$ on which are defined both the MFF $X$ and the Brownian motion
$B$. On this measurable space are defined the probability measures $\Pb
_x=\Pb^X\otimes\Pb^B_x$ (with expectation $\E_x$) for all $x\in\R^2$.
Notice that under $\Pb_x$, the MFF $X$ and the Brownian motion are
independent. We will also denote by $\mathcal{F}_n$ the sigma-algebra
generated by the fields $(Y_k)_{k\leq n}$, that is, $\mathcal
{F}_n=\sigma\{Y_k(x);k\leq n,x\in\R^2\}$. Finally, we mention that we
will sometimes consider other Brownian motions $\bar{B},W$: the
convention of notation will be the same as for $B$.
\end{notation}

\subsection{Gaussian multiplicative chaos}
Let us fix $\gamma\geq0$. We consider the random measure for $n\geq
1$ (the constant $c_n$ is defined in the previous subsection)
%
\begin{equation}
\label{def:Mn} M_n(dx)=c_n^{-\sfrac{\gamma^2}{2}}e^{ \gamma X_n(x)}
\,dx,
\end{equation}
defined on the Borel sets of $\R^2$, which will be called
$n$-regularized Liouville measure. Classical theory of Gaussian
multiplicative chaos (\cite{cfKah} or \cite{review}, Theorem 2.5)
ensures that, $\Pb^X$ almost surely, the family $(M_n)_{n\geq1}$
weakly converges as $n\to\infty$ toward a limiting Radon measure $M$,
which is called the Liouville measure. The limiting measure is
nontrivial if and only if $\gamma\in[0,2)$. We will denote by $\xi_M$
the power law spectrum of $M$ (see \cite{Rnew1,bacry,Rnew10}, e.g.):
%
\begin{equation}
\label{def:xim} \forall p\geq0,\qquad\xi_M(p)=\biggl(2+
\frac{\gamma^2}{2}\biggr)p-\frac
{\gamma^2}{2}p^2.
\end{equation}
Recall (see \cite{cfKah} or \cite{review}, Theorems 2.11 and 2.12)
that for all bounded Borel set $A$ and $p< 4/\gamma^2$, we have $\E
[M(A)^p]<+\infty$ and that
%
\begin{equation}
\label{eq:pls} \sup_{r<1}r^{- \xi_M(p)}\E\bigl[M(r
A)^p\bigr]\leq C_p
\end{equation}
for some constant $C_p$ only depending on $p$.\vadjust{\goodbreak}

Let us emphasize that Kahane's theory of Gaussian multiplicative chaos
ensures that the law of the measure $M$ does not depend on the chosen
regularization $(X_n)_n$ of $X$ (see \cite{cfKah,cfRoVa,review}).
Furthermore, in the case of a GFF $X$ on a planar domain, this result
is reinforced in \cite{cfDuSh} for the Liouville measure: the authors
prove that circle average approximations of $X$ and projections of $X$
along any $H^1$ basis yields almost surely the same Liouville measure.
For more recent results on existence and uniqueness, see \cite{shamov}.

We state below a result about the local modulus of continuity of the
measure $M$ as well as its approximating sequence $(M_n)_n$.

\begin{theorem} \label{app:multform}
We set $\alpha=2(1-\frac{\gamma}{2})^2>0$. Let $\epsilon>0$ and $R>0$.
$\Pb^X$-almost surely, there exists a random constant $C>0$ such that:
\[
\sup_{r\in(0,1)}\sup_{x \in[-R,R]^2} \sup
_{n\geq1}r^{-\alpha
+\epsilon
} \bigl(M_n\bigl( B(x,r) \bigr)\bigr)+M
\bigl( B(x,r) \bigr) \leq C.
\]
%
\end{theorem}

\begin{pf}
It suffices to prove the result only for $M$
because $(M_n)_n$ is a martingale converging a.s. toward $M$ (use
Doob's inequalities to estimate the $\sup_n M_n$ in terms of $M$).

We take $R=\frac{1}{2}$ for simplicity. Now, we partition $[-\frac
{1}{2},\frac{1}{2}]^2$ into $2^{2n}$ dyadic squares $(I_n^j)_{1 \leq j
\leq2^{2n}}$ of equal size. If $p$ belongs to $]0, \frac{4}{\gamma
^2}[$, we get
\begin{eqnarray*}
\Pb^X\biggl( \sup_{1 \leq j \leq2^{2n}} M\bigl(I_n^j
\bigr) \geq\frac
{1}{2^{(\alpha
-\epsilon)n}} \biggr) & \leq&2^{p(\alpha-\epsilon)n} \E^X
\biggl[ \sum_{1 \leq j
\leq2^{2n}} M\bigl(I_n^j
\bigr)^p \biggr]
\\
& \leq&\frac{C_p}{2^{(\xi_M(p)-2- (\alpha-\epsilon) p )n}}.
\end{eqnarray*}
By taking $p= \frac{2}{\gamma}$ in the above inequalities [i.e., $\xi
_M(p)-2- (\alpha-\epsilon) p >0$] and by using Borel--Cantelli's lemma,
we obtain that, $\Pb^X$ almost surely, there exists a random constant
$C$ such that
\[
\sup_{1 \leq j \leq2^{2n}} M\bigl(I_n^j \bigr) \leq
\frac{C}{2^{(\alpha
-\epsilon
)n}} \qquad\forall n \geq1.
\]

We conclude by the fact that each ball $ B(x,r) $ is contained in at
most 4 dyadic squares $(I_{n}^j)_{1 \leq j \leq2^{2n}}$ when we choose
$n$ such that $\frac{1}{2^{n+1}} < r \leq\frac{1}{2^{n}}$.
\end{pf}

\subsection{Potential of the measure $M$}

For each $R>0$, let us introduce the Green function $G_R$ of the
Laplacian on the ball $B(0,R)$ with Dirichlet boundary conditions, that
is, ($\delta_x$ stands for the Dirac mass at $x$)
%
\begin{equation}
\triangle G_R(x,\cdot)=-2 \delta_x(\cdot),\qquad
G_R(x,\cdot )_{\mid\partial B(0,R)}=0.
\end{equation}
Keep in mind the distinction between the massive Green function $G_m$
defined by~(\ref{MFF1}) on $\R^2$ and the Green function $G_R$ on a
ball $B(0,R)$.
Despite the similar notation,
this should bring no confusion as we will always refer to the massive
Green function when the subscript is $m$ and the Green function on
balls when the subscript is $R$.

We introduce the $R$-potential of a Borel measure $\mu$ on $\R^2$ by
\[
\forall x\in B(0,R),\qquad g_R (\mu) (x):= \int
_{B(0,R)} G_R(x,y) \mu (dy).
\]
Let us consider the set of measures
\[
\mathcal{M}=\{M_n;n\geq1\}\cup\{M\}.
\]
Furthermore, for $x\in\R^2$, we denote by $M^z$ the shifted measure
$M^z(\cdot)=M(z+\cdot)$.

Now we use Theorem \ref{app:multform} to prove the following.

\begin{proposition}\label{ass:bass}
For any $R>0$, $\Pb^X$-almost surely, we have:
\begin{longlist}[2.]
\item[1.] $\sup_{\mu\in\mathcal{M}}\sup_{x\in B(0,R)}g_R(\mu
)(x)<+\infty$,

\item[2.]  for any $\mu\in\mathcal{M}$, the mapping $x\in\bar
{B}(0,R)\mapsto g_R (\mu) (x)$ is continuous,

\item[3.] $\sup_{x\in B(0,R)}\llvert  g_R(M_n)(x)-g_R(M)(x)\rrvert  \to0$ as $n\to\infty$,

\item[4.]  for any $z_0\in\R^2$, $\lim_{z\to z_0}\sup_{x\in
B(0,R)}\llvert  g_R(M^z)(x)-g_R(M^{z_0})(x)\rrvert  = 0$.
\end{longlist}
\end{proposition}

\begin{pf}
Recall that the Green function $G_R$ satisfies
for all $x,y\in\bar{B}(0,R)$
%
\begin{equation}
\label{ineg:green} G_R(x,y)\leq\frac{1}{\pi} \ln
\frac{1}{\llvert  x-y\rrvert  }+C
\end{equation}
for some constant $C$. Therefore, it suffices to prove that
\[
\sup_{\mu\in\mathcal{M}}\sup_{x\in B(0,R)} \int
_{B(0,R)} \ln \frac
{1}{\llvert  x-y\rrvert  }\mu(dy)<+\infty.
\]
From Theorem \ref{app:multform}, we can find a constant $C$ and
$\alpha
>0$ (depending on $\gamma$ and $R$) such that for all $x \in\bar
{B}(0,2 R)$ and all $r\in(0,R)$
\[
\sup_{\mu\in\mathcal{M}}\mu\bigl(B(x,r)\bigr)\leq Cr^\alpha.
\]
For $\mu\in\mathcal{M}$, we have
%
\begin{eqnarray}\label{eqagain}
\int_{B(0,R)} \ln\frac{1}{\llvert  x-y\rrvert  }\mu(dy)&\leq&\sum
_{n\geq0}\int_{B(0,R)\cap\{2^{-n}R<\llvert  x-y\rrvert  \leq2^{-n+1}R\}} \ln\frac{1}{\llvert  x-y\rrvert  } \mu
(dy)
\nonumber
\\
&\leq& \sum_{n\geq0} ( n\ln2-\ln
R) \mu\bigl(B\bigl(x,2^{-n+1}R\bigr)\bigr)
\\
&\leq& C2^\alpha\sum_{n\geq1} ( n\ln2-\ln R)
R^\alpha2^{-\alpha
n}.\nonumber
\end{eqnarray}
This latter quantity is finite and does not depend on $x$ or $\mu\in
\mathcal{M}$. This proves the first part of our statement.

For the second statement, consider a function $\theta:\R\to\R$ such
that $0\leq\theta\leq1$, $\theta(x)=1$ for $\llvert  x\rrvert  \leq1$ and $\theta
(x)=0$ for $\llvert  x\rrvert  \geq2$. For $\delta>0$ set $\theta_\delta(x)=\theta
(x/\delta)$ and $\bar{\theta}_\delta(x)=1-\theta_\delta(x)$.
Choose $\mu
\in\mathcal{M}$. Observe that for all $\delta>0$
%
\begin{eqnarray} \label{def:gR}
g_R (\mu) (x)&= &\int_{B(0,R)}
G_R(x,y)\theta_\delta(x-y) \mu (dy)\nonumber
\\
&&{} +\int
_{B(0,R)} G_R(x,y)\bar{\theta}_\delta(x-y)
\mu(dy)
\\
&=:& A_{\mu,\delta}(x)+D_{\mu,\delta}(x).\nonumber
\end{eqnarray}
We will show that the mappings $x\mapsto D_{\mu,\delta}(x)$ converge
uniformly toward $x\mapsto g_R (\mu) (x)$ as $\delta\to0$. As
it is obvious to check that the mapping $x\mapsto D_{\mu,\delta}(x)$
is continuous with the help of standard theorems of continuity for
parameterized integrals, this will show that $x\mapsto g_R (\mu) (x)$
is continuous. So, let us show that the family of mappings $x\mapsto
A_{\mu,\delta}(x)$ converges uniformly toward $0$ on $\bar{B}(0,R)$ as
$\delta\to0$. From~(\ref{ineg:green}) again, it is enough to show that
\[
\sup_{x\in\bar{B}(0,R)} \int_{B(0,R)} \ln
\frac{1}{\llvert  x-y\rrvert  } \theta _\delta (x-y) \mu(dy)\to0\qquad\mbox{as }\delta
\to0.
\]
In fact, we will prove a stronger statement. With computations similar
to~(\ref{eqagain}), we get
%
\begin{eqnarray}\label{eq:again2}
&& \sup_{\mu\in\mathcal{M}}\int_{B(0,R)} \ln\frac{1}{\llvert  x-y\rrvert  }
\theta _\delta (x-y)\mu(dy)
\nonumber
\\
&&\qquad   \leq\sup_{\mu\in\mathcal{M}}\sum
_{n\geq\afrac{\ln
2}{-\ln4\delta}}\int_{B(0,R)\cap\{2^{-n-1}<\llvert  x-y\rrvert  \leq2^{-n}\}} \ln \frac{1}{\llvert  x-y\rrvert  }
\mu(dy)
\\
&&\qquad  \leq C\ln2 \sum_{n\geq\afrac{\ln2}{-\ln4\delta}} (n+1)
2^{-\alpha n}.\nonumber
\end{eqnarray}
This latter series converges to $0$ as $\delta\to0$. The function
$g_R(\mu)$ is thus continuous as a uniform limit of continuous functions.

We now prove the third statement. Sticking to the previous notation
(\ref{def:gR}), we have
\begin{eqnarray*}
g_R (M_n) (x)& =: &A_{M_n,\delta}(x)+D_{M_n,\delta}(x).
\end{eqnarray*}
From~(\ref{eq:again2}), we\vspace*{1pt} have $\sup_{\mu\in\mathcal{M}}\sup_{x\in\bar
{B}(0,R)}\llvert  A_{\mu,\delta}(x)\rrvert  \to0$ as $\delta\to0$. Therefore, it
suffices to prove that for each fixed $\delta>0$, the family
$(D_{M_n,\delta})_n$ converges uniformly on $ \bar{B}(0,R)$ toward
$D_{M,\delta}$. Point-wise convergence is ensured by the weak
convergence of the family of measure $(M_n)_n$ toward $M$ as $n\to
\infty$. We just have to show that the family $(D_{M_n,\delta})_n$ is
relatively compact for the topology of uniform convergence on $ \bar
{B}(0,R)$. For each fixed $\delta>0$, the mapping $(x,y)\mapsto
G_R(x,y)\bar{\theta}_\delta(x-y)$ is continuous on $\bar{B}(0,R)^2$
and, therefore, uniformly continuous. The quantity
\[
\omega(\eta)=\sup_{\llvert  x-x' \rrvert  \leq\eta,(x,x',y)\in\bar
{B}(0,R)^3}\bigl\llvert G_R(x,y)
\bar{\theta}_\delta(x-y)-G_R\bigl(x',y\bigr)
\bar{\theta }_\delta\bigl(x'-y\bigr)\bigr\rrvert
\]
thus converges toward $0$ as $\eta\to0$.
We have for $x,x'\in\bar{B}(0,R)^2$
\[
\bigl\llvert D_{M_n,\delta}(x)-D_{M_n,\delta}\bigl(x'\bigr)
\bigr\rrvert \leq\omega\bigl(\bigl\llvert x-x' \bigr\rrvert
\bigr)M_n\bigl(\bar{B}(0,R)\bigr).
\]
We complete the proof with the Arzel\`a--Ascoli criterion and the
relation $\sup_nM_n(\bar{B}(0,R))<+\infty$ almost surely. The proof of
item 4 can be handled the same way, so we let the reader check the
details.
\end{pf}

\subsection{Approximation and construction of the PCAF of $M$}\label{sub:pcaf}

For $R>0$, we further introduce the stopping time
\[
T_R= \inf\bigl\lbrace t>0; B_t \notin B(0,R) \bigr
\rbrace.
\]
Observe that Theorem \ref{app:multform} (or Proposition \ref{ass:bass})
implies that each $\mu\in\mathcal{M}$ does not charge any polar set.
Following \cite{brosamler}, Proposition 3.2, we deduce that $\Pb^X$
a.s. we can associate to each $\mu\in\mathcal{M}$ a unique PCAF
$(F^{\mu,R}_t)_t$ such that the process
\[
\forall t\geq0,\qquad g_R(\mu) (B_{t\wedge T_R})-g_R(
\mu) (B_0)+F^{\mu,R}_t
\]
is a mean zero martingale under $\Pb^B_x$ for all $x\in B(0,R)$. The
Revuz measure of the PCAF $F^{\mu,R}$ is $\mu$ [restricted to the ball
$B(0,R)$].

\begin{notation}\label{not:fn}
When $\mu=M_n$ for some $n\geq1$, we write $F^{n}_R(t)$ instead of the
heavy notation $F^{M_{n},R}_t$. Similarly, we write $F_R(t)$ for $F^{M,R}_t$.
\end{notation}

It may be worth mentioning that we have the explicit expression
%
\begin{equation}
\label{def:Fnz} F^{n}_R(t)= c_{n}^{-\sfrac{\gamma^2}{2}}
\int_0^{t\wedge T_R}e^{\gamma
X_n(B_r)} \,dr.
\end{equation}

The purpose of what follows is now to establish the convergence of the
family of PCAFs $(F^n_R)_n$ toward $F_R$. Following \cite{BaKo}, we
consider the following distance between two Radon measures $\mu,\nu$ on
$\bar{B}(0,R)$:
\[
d_R(\mu, \nu)= \sup_{x \in\bar{B}(0,R)} \bigl\llvert
g_R(\mu) (x)-g_R(\nu) (x) \bigr\rrvert.
\]

Now we state the following lemma; the proof of which is omitted as a
straightforward adaptation of \cite{BaKo}, Proposition 2.1.

%
\begin{lemma}\label{prop:ko}
For all $R,\eta>0$, $x\in\bar{B}(0,R)$ and $\mu_1,\mu_2$ two Borel
measures such that $d_R(\mu_1,\mu_2)\leq1$, we have
\[
\Pb^B_x\Bigl(\sup_{t \geq0} \bigl\llvert
F^{\mu_1,R}_t- F^{\mu_2,R}_t\bigr\rrvert \geq
\eta\Bigr) \leq c_R \exp \biggl(-\frac{\eta}{c_R\sqrt{d_R(\mu_1,\mu_2)}} \biggr)
\]
for some constant $c_R$ that only depends (increasingly) on
\[
\sup_{i=1,2}\sup_{x\in B(0,R)}g_R(\mu_i)(x).
\]
\end{lemma}

From Proposition \ref{ass:bass} item 3, $\Pb^X$-almost surely and for
all $R>0$, we have $d_R(M_n, M)\to0$ as $n\to\infty$. We deduce the
following.

\begin{corollary}\label{cor:surR}
$\Pb^X$-almost surely, for all $x\in\bar{B}(0,R)$, we have
\[
\Pb^B_x\Bigl( \sup_{t \geq0} \bigl
\llvert F^n_R(t)- F_R(t)\bigr\rrvert \geq
\eta\Bigr) \to0,\qquad \mbox{as }n\to\infty.
\]
\end{corollary}

\begin{theorem}\label{th:PCAF}
$\Pb^X$-almost surely, there exists a unique PCAF denoted by $F$ such that
\[
F(t)=F_R(t),\qquad\mbox{for } t<T_R.
\]
Furthermore, $\Pb^X$-almost surely:
\begin{longlist}[1.]
\item[1.] the Revuz measure of $F$ is $M$,
\item[2.] for all $x\in\R^2$ and $T>0$, $\Pb^B_x ( \sup_{t \leq T}
\llvert   c_{n}^{-\sfrac{\gamma^2}{2}}\int_0^{t}e^{\gamma X_n(B_r)} \,dr-
F(t)\rrvert  \geq\eta ) \to0$ as $n\to\infty$,
\item[3.] for all $x\in\R^2$, $\Pb^B_x$-a.s., $F$ is strictly increasing,\vspace*{1pt}
\item[4.] for all $x\in\R^2$, $\Pb^B_x$-a.s., $\lim_{t\to\infty
}F(t)=+\infty$,
\item[5.] the law of the pair $(B,F)$ under $\Pb^B_x$ on the space of
continuous functions on $\R_+$ equipped with the topology of uniform
convergence over compact sets is a continuous function of $x$, meaning
\[
\lim_{x\to x_0}\E^B_x\bigl[G(B,F)\bigr]=
\E^B_{x_0}\bigl[G(B,F)\bigr]
\]
for every bounded continuous function on $C([0,T],\R_+)$ (for $T>0$)
and $x_0\in\R^2$.
\end{longlist}
\end{theorem}

\begin{pf}
Existence and uniqueness of such a
PCAF is a straightforward consequence of the previous results.

Item 2 results from Corollary \ref{cor:surR} provided that
one takes $R$ large enough to make $\Pb^B_x(T_R<T)$ arbitrarily small.

Item 1 results from \cite{fuku}, Theorem 5.1.3 and Lemma 5.1.10, and \cite{brosamler}, Proposition 3.2. Indeed, this shows
that $F$ coincides with the whole plane PCAF of the measure $M$ (recall
that it does not charge polar sets).

Now we focus on items 3 and 4. Obviously, $\Pb^X$-almost
surely and for all $x\in\R^2$, the mapping $t\in\R_+\mapsto F(t)$ is
increasing $\Pb^B_x$-almost surely. This mapping thus defines a measure
on $\R_+$, which we still denote by $F$ with a slight abuse of notation.

From now on, we will use a few auxiliary lemmas along the main
argument: their proofs are postponed after that of Theorem \ref
{th:PCAF}. Recall the definition of $(\mathcal{F}_n)_n$ in Section~\ref{index}.

\begin{lemma}\label{kodtc1}
For each fixed $x\in\R^2$, $\Pb^B_x$-almost surely, for all $t\geq0$
and \mbox{$R>0$}, the family $(F_R^n(t))_n$ is a uniformly integrable
martingale with respect to the filtration $(\mathcal{F}_n)_n$, which
converges $\Pb^X$-almost surely toward $F_R(t)$. We have, $\Pb
^B_x$-almost surely, $\E^X[F_R(t)]=t\wedge T_R$.
\end{lemma}

We prove item 3. We fix $x\in\R^2$ and we first prove that
$\Pb_x$-a.s. $F$ is strictly increasing. It suffices to prove that it
is strictly increasing on $[0,T_R[$ for all $R>0$. We consider a
nonempty interval $I=[s,t]$ with $t<T_R$. $\Pb^B_x$-a.s., the event
$\{
F_R(I)>0\}$ is an event belonging to the asymptotic sigma-algebra
generated by the random processes $(Y_n)_n$, that is,
\[
\bigl\{F_R(I)>0\bigr\}\in\bigcap_{N\geq1}
\sigma\bigl\{Y_k(x); x\in\R^2,k\geq N\bigr\}.
\]
As the processes $(Y_k)_k$ are independent, we can use the Kolmogorov
0--1 law to deduce that, $\Pb^B_x$-almost surely, the event $\{
F_R(I)>0\}$ has $\Pb^X$-probability $0$ or $1$. From Lemma \ref
{kodtc1}, we have $\Pb^B_x$-almost surely
\[
\E^X\bigl[F_R(I)\bigr]=t-s>0.
\]
Therefore, $\Pb^B_x$-almost surely, on the event $\{t<T_R\}$, the event
$\{F_R(I)>0\}$ has $\Pb^X$-probability $1$. Then we can consider a
countable family $(I_p)_p$ of intervals generating the Borel sigma
algebra on $[0,T_R[$. $\Pb^B_x$-almost surely, we have $F_R(I_p)>0$ for
all $p$ (and all $R>0$). This shows that $\Pb^X$-a.s., $F$ has full
support $\Pb^B_x$-a.s., which equivalently means that the random
mapping $t\mapsto F(t)$ is strictly increasing.

So far, we have only proved that, for each $x\in\R^2$, there is a
measurable set $S_x\subset\Omega$ such that $\Pb^X(S_x)=1$ and on $
S_x$, $F$ is strictly increasing under $\Pb^B_x$. Now we want to show
that there is a measurable set $S\subset\Omega$ such that $\Pb^X(S)=1$
and for all $x\in\R^2$, on $S$, $F$ is strictly increasing under $\Pb
^B_x$. Clearly, we can find a set $S$ such that $\Pb^X(S)=1$ for all
$x\in\Q^2$ and $F$ is increasing on $S$ under $\Pb^B_x$. Now we explain
a coupling procedure that will serve to complete the proof of item 3 as
well as 4 and 5. Let us consider another Brownian motion $W$
independent of $B,X$ (even if it means enlarging the space $\Omega$, we
may assume that $W$ is defined on the same probability space than the
MFF $X$ and the Brownian motion $B$). We denote by $\Pb^{B,W}_{x,y}$
the probability measure $\Pb^B_x\otimes\Pb^W_y$. We state the
following coupling lemma, the proof of which is rather elementary and
thus left to the reader.

\begin{lemma}\label{lem:coupling}
Let us denote by $\tau_1$ the first time at which the first components
of $B$ and $W$ coincide and by $\tau_2$ the first time at which the
second components coincide after $\tau_1$:
\[
\tau_1=\inf\bigl\{u>0;B^{1}_u=W^{1}_u
\bigr\},\qquad\tau_2=\inf\bigl\{u>\tau _1;B^{2}_u=W^{2}_u
\bigr\}.
\]
Under $\Pb^{B,W}_{x,y}$, the random process $\bar{B}$ defined by
\[
\bar{B}_t=\cases{ \bigl(W^{1}_t,W^{2}_t
\bigr), &\quad if $t\leq\tau_1$,
\vspace*{3pt}\cr
\bigl(B^{1}_{t},
W^{2}_t\bigr), &\quad if $\tau_1<t\leq
\tau_2$,
\vspace*{3pt}\cr
\bigl(B^{1}_t,B^{2}_t
\bigr), &\quad if $\tau_2<t$,}
\]
is a Brownian motion on $\R^2$ starting from $y$, and coincides with
$W$ for all times $t>\tau_2$. Furthermore, we have
\begin{eqnarray*}
\forall\eta>0,\qquad \lim_{\delta\to0} \sup_{x,y\in\R^2;\llvert  x-y\rrvert  \leq
\delta}\Pb
^{B,X}_{x,y}(\tau_2>\eta)\to0\quad\mbox{and}\quad
\Pb^{B,X}_{x,y}(\tau_2<\infty)=1.
\end{eqnarray*}
\end{lemma}

We can associate $\Pb^X$-a.s. to the Brownian motion $\bar{B}$ a
PCAF, denoted by $F(\bar{B},t)$ to distinguish it from $F$, with
Revuz measure $M$ as prescribed in the beginning of the proof of
Theorem \ref{th:PCAF}. It is also plain to check that $\Pb^X$-a.s., for
all $x,y\in\R^2$, under $\Pb^{B,W}_{x,y}$, the marginal laws of $(B,F)$
and $(\bar{B},F(\bar{B},\cdot))$, respectively, coincide with
the law of $(B,F)$ under $\Pb^B_x$ and $\Pb^B_y$.

Therefore, on $S$ and for $y\in\Q^2$, $\Pb^{B,W}_{x,y}$ a.s., the PCAF
$F(\bar{B},\cdot)$ is strictly increasing on $\R_+$. Furthermore,
the coupling procedure (Lemma \ref{lem:coupling}) also entails that the
mappings $s\in[\tau_2,+\infty[\,\mapsto F(s)-F(\tau_2)$ and $s\in
[\tau
_2,+\infty[\,\mapsto F(\bar{B},s)-F(\bar{B},\tau_2)$ are equal
$\Pb^{B,W}_{x,y}$ a.s. Therefore,\vspace*{1pt} for any $x\in\R^2$ and $y\in\Q^2$,
the above discussion shows that $F(s)$ is strictly increasing for
$s>\tau_2$. If $y$ is chosen arbitrarily close to $x$, Lemma \ref
{lem:coupling} shows that $\tau_2\to0$ in probability. We deduce that
$F$ is strictly increasing on $\R_+$. This completes the proof of item
3 as we have shown that on $S$ for all $x\in\R^2$, $F$ is strictly
increasing $\Pb^{B}_{x}$ a.s.

We now prove item 4. Once again, the coupling procedure
(Lemma \ref{lem:coupling}) shows that it is enough to prove item 4
$\Pb
^X$ a.s. for only one $x\in\R^2$: because $\Pb^{B,W}_{x,y}(\tau
_2<\infty
)=1$ for all $x,y$, it is plain to deduce that $\Pb^B_x(F(t)\to\infty
\mbox{ as }t\to\infty)=\Pb^B_y(F(t)\to\infty\mbox{ as
}t\to\infty
)$ for all $y\in\R^2$.
So we work under $\Pb^B_0$.

We consider the following sequence of stopping times associated to the
Brownian motion:
\[
T_n = \inf\bigl\lbrace t> 0, \llvert B_t\rrvert = 2n
\bigr\rbrace,\qquad\bar{T}_n= \inf \bigl\lbrace t> T_n,
\llvert B_t-B_{T_n} \rrvert = \tfrac{1}{4} \bigr
\rbrace.
\]
We also consider an increasing sequence of integers $(n_j)_{j \geq1}$
such that the following property holds for all $l \leq k$:
%
\begin{eqnarray}\label{loveko}
\sum_{l \leq j<j' \leq k} \alpha_{j,j'}
\leq k-l+1,
\nonumber\\[-8pt]\\[-8pt]
\eqntext{\displaystyle \mbox{with } \alpha_{j,j'}= \sup_{\llvert  x\rrvert   \leq2n_j+1/4, \llvert  y\rrvert   \geq2 n_{j'}-1/4}
G_m(x,y).}
\end{eqnarray}
Such a sequence exists because $G_m$ defined by~(\ref{MFF1}) satisfies
$G_m(x,y)\leq c e^{-c\llvert  x-y\rrvert  }$ for $\llvert  x-y\rrvert  \geq1$ and some constant $c>0$.
Recall that we identify $F$ and its associated measure. From the Markov
inequality and Fatou's lemma, we obtain
%
\begin{eqnarray}\label{expect}
&& \Pb_0 \biggl( \bigcap_{l \leq j \leq k } \bigl\{
F\bigl(]T_{n_j}, \bar{T}_{n_j}]\bigr) \leq c\bigr\} \biggr) \nonumber
\\
&&\qquad \leq
c^{k-l+1} \E_0 \biggl[ \prod_{l \leq j \leq k}
\bigl(F\bigl(]T_{n_j}, \bar{T}_{n_j}]\bigr)
\bigr)^{-1} \biggr]
\nonumber\\[-8pt]\\[-8pt]\nonumber
&&\qquad  \leq \liminf_{n\to\infty} c^{k-l+1}
\E_0 \biggl[ \biggl( \int_{ \prod_{l \leq j \leq k} ]T_{n_j}, \bar{T}_{n_j}] }
c_n^{-(k-l+1)\sklfrac{\gamma^2}{2}}
\\
&&\quad\qquad{}\times
e^{ \gamma(X_n(B_{s_l}) + \cdots+
X_n(B_{s_k}) ) }\,ds_l\cdots
ds_k \biggr)^{-1} \biggr].\nonumber 
\end{eqnarray}
We want to get rid of the long range correlations of the MFF $X$. To
this purpose, we introduce a log-correlated random distribution $\bar
{X}$ with covariance kernel $\E[\bar{X}(x)\bar{X}(y)]=\ln_+\frac
{1}{\llvert  y-x \rrvert  }$. It is a kernel of $\sigma$-positive type \cite{review}, Proposition~2.15. We can find an\vspace*{1pt} approximation family ``\`a la Kahane''
$(\bar{X}_n)_n$ of $\bar{X}$ such that $\E[\bar{X}_n(x)\bar
{X}_n(y)]\leq\ln_+\frac{1}{\llvert  y-x \rrvert  }$ and
%
\begin{equation}
\E\bigl[\bar{X}_n(x)\bar{X}_n(y)\bigr]-D \leq\E
\bigl[X_n(x)X_n(y)\bigr]\leq\E\bigl[\bar
{X}_n(x)\bar{X}_n(y)\bigr]+D
\end{equation}
for some constant $D$ which does not depend on relevant quantities (in
particular not on $n$, e.g., \cite{review}, proof of Proposition~2.15). By using this relation and~(\ref{loveko}), we get for
all $(x_l,\dots,x_k),(y_l,\dots,y_k)\in\prod_{l \leq j \leq k}
]T_{n_j}, \bar{T}_{n_j}] $
%
\begin{eqnarray}\label{porcasse}
&& \E^{X} \Biggl[ \Biggl( \sum_{j=l}^kX_n(x_j)
\Biggr) \Biggl(\sum_{j=l}^kX_n(y_j)
\Biggr) \Biggr]
\nonumber
\\
&&\qquad   =\sum_{j=l}^k
\E^{X} \bigl[ X_n(x_j)X_n(y_{j})
\bigr]+\sum_{j,j'=l,j\neq j'}^k\E^{X}
\bigl[ X_n(x_j)X_n(y_{j'}) \bigr]
\nonumber\\[-8pt]\\[-8pt]\nonumber
&&\qquad   \leq(k-l+1)D+\sum_{j=l}^k
\E^{X} \bigl[ \bar {X}_n(x_j)\bar
{X}_n(y_{j}) \bigr]+2\sum_{l \leq j<j' \leq k}
\alpha_{j,j'}
\nonumber
\\
&&\qquad   \leq(k-l+1) (D+2)+\E^{\bar{X}} \Biggl[ \Biggl(\sum
_{j=l}^k\bar {X}_n(x_j)
\Biggr) \Biggl(\sum_{j=l}^k
\bar{X}_n(y_j) \Biggr) \Biggr].\nonumber
\end{eqnarray}
In the last line, we have used the fact that $\E^{\bar{X}} [ \bar
{X}_n(x_j)\bar{X}_n(y_{j'}) ]=0$ if $j\neq j'$.

Now we want to apply Lemma \ref{cvx} with $Y(x_l,\dots,x_k)=
\sum_{j=l}^kX_n(x_j)$ (with kernel $K$) and $Y'(x_l,\break \dots,x_k)=\sum_{j=l}^k\bar{X}_n(x_j)$ (with kernel $K'$),  $\nu(dx_l,
\dots,dx_k)=\mu
_l(dx_l)\times\cdots\times\mu_k(dx_k)$ where each $\mu_l$ stands
for the
occupation measure of the Brownian motion between the times $T_{n_j}$
and $\bar{T}_{n_j}$ and the convex function $x\mapsto1/x$ (in fact, this
function is discontinuous at $0$ but this is not a problem: truncate it
in order to have a continuous convex function, apply Kahane's
inequality and then remove the truncation). We have shown above that
$K\leq K'+C$ with $C=(k-l+1)(D+2)$. The last point is that the
exponential term in the expectation~(\ref{expect}) is not renormalized
by the variance. However, the above inequality shows that
\[
\E^{X} \Biggl[ \Biggl(\sum_{j=l}^kX_n(x_j)
\Biggr)^2 \Biggr]\leq\sum_{j=l}^k
\E ^X\bigl[X_n(x_j)^2
\bigr]+2D(k-l+1)
\]
in such a way that $(k-l+1) \ln c_n\leq\E^{X} [ (\sum_{j=l}^kX_n(x_j) )^2 ]\leq(k-l+1) \ln c_n+2D(k-l+1)$.

Hence, even if it means multiplying the constant $c$ by a deterministic
constant that does not depend on $k,l$, we can replace the term
$c_n^{-(k-l+1)\sklfrac{\gamma^2}{2}}$ in~(\ref{expect}) by $\exp
(-\frac
{\gamma^2}{2}\E^{X} [ (\sum_{j=l}^kX_n(x_j) )^2
] )$. We
are then in position to apply Lemma \ref{cvx} item~2, which tells us
that we can replace $X_n$ in~(\ref{porcasse}) by $\bar{X}_n$ at the
cost of replacing the constant $c$ by another constant $c'$, which
still does not depend on $n$ (only on $D,\gamma$). The main advantage
of this procedure is that we deal now with a field $\bar{X}_n$ that
possesses strong decorrelation properties: if, for a set $A\subset\R
^2$, we denote by $\mathcal{F}_{n,A}$ the sigma algebra generated by
the random variables $\{\bar{X}_n(x), x\in A\}$ then $\mathcal
{F}_{n,A}$ is independent of $\mathcal{F}_{n,B}$ as soon as $\operatorname{dist}(A,B)>1$.

In what follows, we still stick to the notation $\E_0$ to denote
expectation with respect to the probability measure $\Pb^{\bar
{X}}\otimes\Pb^B_0$ and $\bar{c}_n=\E^{\bar{X}}[\bar{X}_n(x)^2]$,
which does not depend on $x$ by stationarity. We have by using the
strong Markov property of the Brownian motion and the fact that $\bar
{X}_n$ is decorrelated at distance $1$
\begin{eqnarray*}
&& \Pb_0 \biggl( \bigcap_{l \leq j \leq k } \bigl\{
F\bigl(]T_{n_j}, \bar{T}_{n_j}]\bigr) \leq c\bigr\} \biggr)
\\
&&\qquad \leq \liminf_{n\to\infty}
\bigl(c'\bigr)^{k-l+1} \E_0 \biggl[ \biggl(
\int_{
\prod
_{l \leq j \leq k} ]T_{n_j}, \bar{T}_{n_j}] } \bar {c}_n^{-(k-l+1)\sklfrac
{\gamma^2}{2}}
\\
&&\quad\qquad{}\times e^{ \gamma(\bar{X}_n(B_{s_l}) + \cdots+ \bar
{X}_n(B_{s_k}) ) }
\,ds_l\cdots ds_k \biggr)^{-1} \biggr]
\\
&&\qquad =  \liminf_{n\to\infty} \bigl(c'
\bigr)^{k-l+1} \E_0 \biggl[ \biggl( \int_{
]T_{n_1}, \bar{T}_{n_1}]}
\bar{c}_n^{- \sfrac{\gamma^2}{2}} e^{
\gamma
\bar{X}_n(B_{s}) } \,ds \biggr)^{-1}
\biggr]^{k-l+1}.
\end{eqnarray*}
Notice\vspace*{1pt} that $n\mapsto\int_{ ]T_{n_1}, \bar{T}_{n_1}]} \bar{c}_n^{- \sfrac{\gamma^2}{2}} e^{ \gamma\bar{X}_n(B_{s}) } \,ds $ is a uniformly
integrable martingale and converges toward a random variable denoted by
\[
\int_{ ]T_{n_1}, \bar{T}_{n_1}]} e^{ \gamma\bar{X} (B_{s})
-\sklfrac
{\gamma^2}{2} \E_0[\bar{X}^2(B_s)] } \,ds
\]
(the proof is identical to
Lemma \ref{kodtc1} and the notation is due to the fact that the limit
is a Gaussian multiplicative chaos). Therefore, the Jensen inequality
leads to
\begin{eqnarray*}
&& \Pb_0 \biggl( \bigcap_{l \leq j \leq k } \bigl\{
F\bigl(]T_{n_j}, \bar{T}_{n_j}]\bigr) \leq c\bigr\} \biggr)
\\
&&\qquad \leq
\biggl( c' \E_0 \biggl[ \biggl( \int
_{ ]T_{n_1}, \bar{T}_{n_1}]} e^{
\gamma\bar{X} (B_{s}) -\sklfrac{\gamma^2}{2} \E_0[\bar{X}^2(B_s)] } \,ds \biggr)^{-1} \biggr]
\biggr)^{k-l+1}.
\end{eqnarray*}
Let us admit for a while that the above expectation in the right-hand
side is finite (this will be proved below in Lemma \ref{negative}).
This inequality shows that we can choose $c$ small enough such that
\[
\Pb_0 \biggl( \bigcap_{l \leq j < \infty} \bigl\{
F\bigl(]T_{n_j}, \bar{T}_{n_j}]\bigr) \leq c\bigr\} \biggr) = 0.
\]
Thus, we get
\[
\Pb_0 \biggl( \bigcap_{l\geq1}\bigcup
_{l \leq j < \infty} \bigl\{ F\bigl(]T_{n_j},
\bar{T}_{n_j}]\bigr) > c\bigr\} \biggr) =1.
\]
Since $\lim_{t\to\infty}F( t) \geq c \sum_{j \geq1} 1_{ \{F(]T_{n_j},
\bar{T}_{n_j}]) > c \}}$,
the proof of item 4 is complete.

It remains to prove item 5. For $y\in\R^2$, we denote by
$M^y$ the shifted measure $M^y(A)=M(A+y)$, and $F^y$ its associated
PCAF. We claim that it is enough to prove that, $\Pb^X$-a.s. for all
$x\in\R^2$ and $\eta>0$
%
\begin{equation}
\label{couplingF} \lim_{y\to x}\Pb^{B}_{0}
\Bigl(\sup_{t\leq T}\bigl\llvert F^y(t)-F^x(t)
\bigr\rrvert \geq\eta \Bigr)=0.
\end{equation}
Indeed, if~(\ref{couplingF}) is true, then for any uniformly continuous
function $G$ (with modulus $m$ and bounded by $K$) on $C(\R_+; \R
^2\times\R_+)$, we have for all $\eta>0$
\begin{eqnarray*}
&& \bigl\llvert \E^B_x\bigl[G(B,F)\bigr]-
\E^B_y\bigl[G(B,F)\bigr]\bigr\rrvert
\\
&&\qquad   =\bigl\llvert \E^B_0\bigl[G
\bigl(x+B,F^x\bigr)\bigr]- \E^B_0\bigl[G
\bigl(y+B,F^y\bigr)\bigr]\bigr\rrvert
\\
&&\qquad   \leq K\ind_{\{\llvert  x-y\rrvert  >\eta\}}+K\Pb^{B}_{0}
\Bigl(\sup_{t\leq
T}\bigl\llvert F^x(t)-F^y(t)
\bigr\rrvert >\eta\Bigr)+m(\eta).
\end{eqnarray*}

We can then pass to the limit as $y\to x$ and use~(\ref{couplingF}) to
prove that the first two terms go to $0$ and then choose $\eta$
arbitrarily small to conclude.

To establish~(\ref{couplingF}), for all $R>0$, we use Lemma \ref
{prop:ko} to get
\begin{eqnarray*}
&& \Pb^{B}_{0} \Bigl(\sup_{t\leq T}\bigl
\llvert F^y(t)-F^x(t)\bigr\rrvert \geq\eta \Bigr)
\\
&&\qquad   =\Pb^{B}_{0} \Bigl(\sup
_{t\leq0}\bigl\llvert F^y(t\wedge
T_R)-F^x(t\wedge T_R)\bigr\rrvert \geq
\eta \Bigr)+\Pb^{B}_{0} (T_R\leq T)
\\
&&\qquad  \leq c_R\exp \biggl(-\frac{\eta}{c_R\sqrt
{d_R(M^y,M^x)}} \biggr)+
\ind_{\{d_R(M^y,M^x)\geq1\}}+\Pb^{B}_{0} (T_R\leq T).
\end{eqnarray*}
Since $d_R( M^y,M^x)\to0$ as $y\to x$ (cf. item 4 of Proposition \ref
{ass:bass}), we deduce
\[
\limsup_{y\to x}\Pb^{B}_{0} \Bigl(\sup
_{t\leq T}\bigl\llvert F^y(t)-F^x(t)\bigr
\rrvert \geq \eta \Bigr)\leq\Pb^{B}_{0} (T_R
\leq T).
\]
We complete the proof by letting $R\to\infty$.
\end{pf}

\begin{pf*}{Proof of Lemma \ref{kodtc1}} $\Pb^B_x$-a.s., the family
$(F^n_R(t))_n$ is a nonnegative martingale w.r.t. the filtration
$(\mathcal{F}_n)_n$ and, therefore, converges almost surely as $n\to
\infty$.

Let us prove that it is uniformly integrable. Denote by $\nu$ the
occupation measure of the Brownian motion $B$ between $0$ and $t\wedge
T_R$. Observe that
\[
\int_{\R^2}c_{n}^{-\sfrac{\gamma^2}{2}}e^{\gamma X_n(z)}
\nu(dz)=F_R^n(t).
\]
From \cite{cfKah} (see also \cite{review}), we just have to prove that
$\Pb^B_x$-a.s.
\[
\int_{\R^2}\int_{\R^2} \frac{1}{\llvert  z-z'\rrvert  ^\alpha}
\nu(dz) \nu \bigl(dz'\bigr)<+\infty
\]
for some $\alpha<2$. This statement is elementary (just compute the
expectation), and thus left to the reader (much stronger statements are
discussed in \cite{dembo}, Section~10, e.g.).
\end{pf*}

\subsection{Study of the moments and power law spectrum} \label{sub.moments}
In this section, we investigate the finiteness of the moments of the
PCAF. We will say that $F$ possesses moments of order $q$ if we have
$\E
_x[F(t)^q]<+\infty$ for all $t>0$ and $x\in\R^2$.

\begin{theorem}[(Positive moments and power law spectrum)]\label{lemmamoments}
(1)~If $\gamma<2$, the mapping $F$ possesses moments of order $q$ for
$0\leq q< 4/\gamma^2$.
(2) If $F$ admits moments of order $q\geq1$ then, for all $s\in[0,1]$
and $t\in[0,T]$:
\[
\E_x \bigl[\bigl(F(t+s)-F(t)\bigr)^q\bigr]\leq
C_q s^{\xi(q)},
\]
where
\[
\xi (q)= \biggl(1+
\frac{\gamma^2}{4} \biggr)q-\frac{\gamma^2}{4}q^2
\]
and $C_q>0$ is some constant independent of $x,T$.
\end{theorem}


\begin{pf} There is here no exception to the rule in
multiplicative chaos theory that studying finiteness of the moments is
technically heavy. So, the entire Appendix~\ref{app:moments} is devoted
to the proof of item 1.


Now we assume that $\gamma<2$ and that $F$ possesses moments of order
$q\geq1$. We prove the estimate concerning the power law spectrum. By
stationarity, we may assume that $x=0$. We first prove it when $t=0$
and then we deduce the uniform estimate in $t$. Under $\Pb_0$, we have
%
\begin{eqnarray}\label{starmom}
F(s)&=&\int_0^{s}e^{\gamma X (B_r)-\sklfrac{\gamma^2}{2}\E[X (B_r)^2]} \,dr
\nonumber
\\
& =& s\int
_0^{1}e^{\gamma X (B_{u s})-\sklfrac{\gamma^2}{2}\E[X (B_{u
s})^2]} \,du
\\
&\stackrel{\mathrm{law}} {=} & s\int_0^{1}e^{\gamma X (\sqrt{s}B_{u})-\sklfrac
{\gamma
^2}{2}\E[X (\sqrt{s}B_{u})^2]}
\,du.\nonumber
\end{eqnarray}
Let us stress here that the above computations are of course only
formal as all the quantities are understood as limits: yet, the final
statement is correct as can be seen by applying the argument to the
same regularized quantities and by passing to the limit to remove the
cutoff. Then, from~(\ref{ssi}), we have $G_m(\sqrt{s}u,\sqrt
{s}v)\leq
\ln\frac{1}{\sqrt{s}}+G_m(u,v)$ for all $u,v\in\R^2$.
Then, by taking the $q$th power and expectation in~(\ref{starmom}) and
Kahane's convexity inequalities (see\vspace*{1pt} Lemma \ref{cvx} item 2), we get
for some Gaussian random variable $\Omega_s$ with mean $0$ and variance
$-\frac{1}{2}\ln s$ and independent of $ \int_0^{t}e^{\gamma X (
B_{u})-\sklfrac{\gamma^2}{2}\E[X ( B_{u})^2]} \,du $
\begin{eqnarray*}
\E_0 \bigl[F(s)^q\bigr]&\leq& s^q
\E_0 \biggl[ \biggl(e^{\gamma\Omega_{s}+\sklfrac
{\gamma
^2}{4}\ln s}\int_0^{1}e^{\gamma X ( B_{u})-\sklfrac{\gamma^2}{2}\E[X (
B_{u})^2]}
\,du \biggr)^q \biggr]
\\
&=&s^q \E_0 \bigl[e^{q \gamma\Omega_{s}+q\sklfrac{\gamma^2}{4}\ln
s} \bigr]\E
_0 \biggl[ \biggl(\int_0^{1}e^{\gamma X ( B_{u})-\sklfrac{\gamma^2}{2}\E[X (
B_{u})^2]}
\,du \biggr)^q \biggr]
\\
&=&C_q s^{\xi(q)},
\end{eqnarray*}
where $C_q=\E_0 [ (\int_0^{1}e^{\gamma X ( B_{u})-\sklfrac
{\gamma
^2}{2}\E[X ( B_{u})^2]} \,du )^q ]$ is independent of $s,x$.

Now we treat the general case $t\neq0$. By using the Markov property
of the Brownian motion and the stationarity of $X$, it is readily seen that
\begin{eqnarray*}
\E_x \bigl[\bigl(F(t+s)-F(t)\bigr)^q\bigr]&=&
\E_0 \bigl[F(s)^q\bigr].
\end{eqnarray*}\upqed
%
\end{pf}

\begin{corollary}
Set $\alpha=(1-\frac{\gamma}{2})^2$. For each $T>0$ and $\epsilon>0$,
there exists a random constant $C>0$ such that $\Pb_x$ a.s.
\[
\sup_{0\leq s<t\leq T}\bigl\llvert F(t)-F(s)\bigr\rrvert \leq C\llvert
t-s\rrvert ^{\alpha-\epsilon}.
\]
\end{corollary}

\begin{pf} Similar to Theorem \ref{app:multform} and thus
left to the reader.\end{pf}

Now we investigate finiteness of moments of negative
order. Denote by $T^x_r$ the first exit time of the Brownian motion $B$
out of the disk $B(x,r)$ for $r\in\,]0,1]$.

\begin{proposition}\label{negative}
For all $q>0$, there exists some constant $C_q>0$ (depending on $q$)
such that for all $x\in\R^2$ and $r\in[0,1]$
%
\begin{equation}
\sup_{n\geq0} \E_x \biggl[ \biggl(
c_n^{- \sfrac{\gamma^2} 2} \int_0^{T^x_r}
e^{ \gamma X_n( B_s)} \,ds \biggr)^{-q} \biggr] \leq C_q
r^{2\xi(-q)}.
\end{equation}
\end{proposition}

\begin{pf} Without loss of generality, we can take $x=0$ by
stationarity of the field $X$. Furthermore, from Kahane's convexity
inequalities Lemma \ref{cvx}, it suffices to prove the result for one
log-correlated Gaussian field with a kernel of $\sigma$-positive type.
Let us choose the exact scale invariant field $\bar{X}$ with covariance
kernel given by
\[
\E\bigl[\bar{X}(x)\bar{X}(y)\bigr]=\ln_+\frac{2}{ \llvert  x-y\rrvert  }
\]
with white noise decomposition $(\bar{X}_\epsilon)_{\epsilon\in\,]0,1]}$
of $\bar{X}$ as constructed in \cite{cfRoVa}. More precisely, the
correlation structure of $(\bar{X}_\epsilon)_{\epsilon\in\,]0,1]}$ is
given for $\epsilon,\epsilon' \in\,]0,1]$ by
\begin{eqnarray*}
\hspace*{-3pt}&& \E\bigl[\bar{X}_\epsilon(x)\bar{X}_{\epsilon'}(y)\bigr]
\\
\hspace*{-3pt}&&\qquad = \cases{ 0, &
\quad if $\llvert x-y\rrvert >2$,
\cr
\displaystyle \ln\frac{2}{ \llvert  x-y\rrvert  }, &\quad if $\max\bigl(
\epsilon,\epsilon'\bigr)\leq \llvert x-y\rrvert \leq 2$,
\cr
\displaystyle \ln
\frac{2}{\max(\epsilon,\epsilon')}+2\biggl( 1- \frac
{\llvert  x-y\rrvert  ^{1/2}}{\max
(\epsilon,\epsilon')^{1/2}} \biggr), &\quad if $\llvert y-x
\rrvert \leq\max\bigl(\epsilon,\epsilon'\bigr)$.}
\end{eqnarray*}
This covariance structure entails some interesting properties that we
detail now. The process $\epsilon\rightarrow\bar{X}_\epsilon$ has
independent increments, meaning that for $\epsilon'>\epsilon$ the field
$\bar{X}_{\epsilon}-\bar{X}_{\epsilon'}$ is independent of the sigma
algebra $\sigma\{X_u(z);z\in\R^2,u\geq\epsilon'\}$. For $\epsilon
'>\epsilon$, the field $\bar{X}_{\epsilon,\epsilon'}:=\bar
{X}_{\epsilon
}-\bar{X}_{\epsilon'}$ has a correlation cutoff of length $\epsilon'$,
meaning that the fields $(\bar{X}_{\epsilon,\epsilon'}(x))_{x \in A}$
and $(\bar{X}_{\epsilon,\epsilon'}(x))_{x \in B}$ are independent
whenever the Euclidean distance between the two sets $A,B$ is greater
than $\epsilon'$. Finally, we have the following relation in law for
$r,\epsilon\leq1$:
%
\begin{equation}
\label{scalingko} \bigl(X_{r \epsilon}(r x)\bigr)_{\llvert  x\rrvert  \leq1}\stackrel{\mathrm{law}}
{=} \bigl(\Omega _r+X_{\epsilon
}( x)\bigr)_{\llvert  x\rrvert  \leq1},
\end{equation}
where $\Omega_{r}$ is a Gaussian distribution $\mathcal{N}(0, \ln
\frac
{1}{r})$ and is independent of $(\bar{X}_\epsilon)_\epsilon$.
In what follows, we stick to the Notation \ref{product} with $\bar{X}$
instead of $X$. Finally, we just write $T_r$ for $T_r^0$.

So, we have to prove
%
\begin{equation}
\label{equintermediaire} \sup_{\epsilon\in\,]0,1]} \E_0 \biggl[ \biggl(
2^{-\sfrac{\gamma^2}
2}\epsilon^{\sfrac{\gamma^2} 2}\int_0^{T_r}
e^{ \gamma\bar{X}_\epsilon( B_s)} \,ds \biggr)^{-q} \biggr] \leq C_q
r^{2\xi(-q)}.
\end{equation}
Notice that the supremum is reached for $\epsilon\to0$ by the
martingale property and the Jensen inequality. Now, if $\widetilde
{T}_{\sfrac {1}{4}}$ is the first time the Brownian motion
$(B_{t+T_{3/4}}-B_{T_{3/4}})_{t\geq0}$ hits the disk of radius $\frac
{1}{4}$, we get for $\epsilon<1/4$
\begin{eqnarray*}
&& 2^{-\sfrac{\gamma^2} 2}\epsilon^{\sfrac{\gamma^2} 2} \int_0^{T_1}
e^{\gamma\bar{X}_{\eps}( B_s)} \,ds
\\
&&\qquad \geq2^{-\sfrac{\gamma^2}
2}\epsilon ^{\sfrac{\gamma^2} 2} \int
_0^{T_{\sfrac{1}{4}}} e^{ \gamma\bar
{X}_{\eps
}( B_s)} \,ds + 2^{-\sfrac{\gamma^2} 2}
\epsilon^{\sfrac{\gamma^2}
2}\int_{0}^{\widetilde{T}_{\sfrac{1}{4}}}
e^{ \gamma\bar{X}_{\eps}(
B_{s+T_{3/4}})} \,ds
\\
&&\qquad \geq8^{-\sfrac{\gamma^2}{2} }e^{\gamma\inf_{\llvert  x\rrvert  \leq1} \bar
{X}_{1/4}(x) } \biggl( (4\epsilon)^{ \sfrac{\gamma^2} 2 }\int
_0^{T_{\sfrac{1}{4}}} e^{ \gamma\bar{X}_{\eps,1/4}( B_s) } \,ds
\\
&&\quad\qquad{}+ (4
\epsilon)^{ \sfrac{\gamma^2} 2 }\int_{0}^{\widetilde{T}_{\sfrac{1}{4}}}
e^{ \gamma\bar{X}_{\eps,1/4}( B_{s+T_{3/4}}) } \,ds \biggr).
\end{eqnarray*}
The main observation is that, under the annealed measure $\E_0$, the
above two integrals are independent identically distributed random
variables. Indeed, by considering two bounded continuous functionals
$F,G$, we get
\begin{eqnarray*}
P(F,G)&:=& \E_0 \biggl[ F \biggl( (4\epsilon)^{ \sfrac{\gamma^2} 2 } \int
_0^{T_{\sfrac{1}{4}}} e^{ \gamma\bar{X}_{\eps,1/4}( B_s)} \,ds \biggr)
\\
&&{}\times G \biggl(
(4\epsilon)^{ \sfrac{\gamma^2} 2 } \int_{0}^{\widetilde
{T}_{\sfrac
{1}{4}}}
e^{ \gamma\bar{X}_{\eps,1/4}( B_{s+T_{3/4}}) } \,ds \biggr) \biggr]
\\
&=& \E^{B}_0 \biggl[ \E^{\bar{X}} \biggl[ F \biggl(
(4\epsilon)^{ \sfrac {\gamma^2} 2 } \int_0^{T_{\sfrac{1}{4}}}
e^{ \gamma\bar{X}_{\eps,1/4}(
B_s) } \,ds \biggr) \biggr]
\\
&&{}\times  \E^{\bar{X}} \biggl[ G \biggl( (4\epsilon)^{ \sfrac{\gamma^2} 2 } \int
_{0}^{\widetilde{T}_{\sfrac{1}{4}}} e^{ \gamma\bar{X}_{\eps,1/4}(
B_{s+T_{3/4}}) } \,ds \biggr) \biggr]
\biggr],
\end{eqnarray*}
where\vspace*{2pt} we have used the fact that $\bar{X}_{\eps,1/4}$ has a correlation
cutoff of length $1/4$. If we use now the stationarity of the field
$\bar{X}_{\eps,1/4}$ and the independence of the increments of the
standard Brownian motion, we obtain
\begin{eqnarray*}
P(F,G) &=& \E^{B}_0 \biggl[ \E^{\bar{X}} \biggl[ F
\biggl( (4\epsilon )^{ \sfrac {\gamma^2} 2 } \int_0^{T_{\sfrac{1}{4}}}
e^{ \gamma\bar{X}_{\eps,1/4}(
B_s) } \,ds \biggr) \biggr]
\\
&&{}\times  \E^{\bar{X}} \biggl[ G \biggl( (4\epsilon)^{ \sfrac{\gamma^2} 2 } \int
_{0}^{\widetilde{T}_{\sfrac{1}{4}}} e^{ \gamma\bar{X}_{\eps,1/4}(
B_{s+T_{3/4}}-B_{T_{3/4}}) } \,ds \biggr) \biggr]
\biggr]
\\
&= &\E_0 \biggl[ F \biggl( (4\epsilon)^{ \sfrac{\gamma^2} 2 } \int
_0^{T_{\sfrac{1}{4}}} e^{ \gamma\bar{X}_{\eps,1/4}( B_s) } \,ds \biggr) \biggr]
\\
&&{}\times  \E_0 \biggl[ G \biggl( (4\epsilon)^{ \sfrac{\gamma^2} 2 } \int
_{0}^{\widetilde{T}_{\sfrac{1}{4}}} e^{ \gamma\bar{X}_{\eps,1/4}(
B_{s+T_{3/4}}-B_{T_{3/4}}) } \,ds \biggr) \biggr]
\\
&=&\E_0 \biggl[ F \biggl( (4\epsilon)^{ \sfrac{\gamma^2} 2 } \int
_0^{T_{\sfrac{1}{4}}} e^{ \gamma\bar{X}_{\eps,1/4}( B_s) } \,ds \biggr) \biggr]
\\
&&{}\times
\E_0 \biggl[ G \biggl( (4\epsilon)^{ \sfrac{\gamma^2} 2 } \int
_0^{T_{\sfrac{1}{4}}} e^{ \gamma\bar{X}_{\eps,1/4}( B_s) } \,ds \biggr) \biggr].
\end{eqnarray*}
Furthermore, for all $r \in\,]0,1]$, we have
\begin{eqnarray*}
(r \eps/2)^{\sfrac{\gamma^2} 2}\int_0^{T_r}
e^{ \gamma\bar{X}_{r
\eps
}(B_s) } \,ds & =& (r \eps/2)^{\sfrac{\gamma^2} 2} r^2 \int
_0^{\sfrac {T_r}{r^2}} e^{ \gamma\bar{X}_{r \eps}(r \sklfrac{B_{r^2 s'}}{r}) } \,ds'
\\
& =& (r \eps/2)^{\sfrac{\gamma^2} 2} r^2 \int_0^{\widetilde
{T}_1}
e^{
\gamma\bar{X}_{r \eps}(r \widetilde{B}_{s'}) } \,ds',
\end{eqnarray*}
where $\widetilde{B}_{s'} = r^{-1}B_{r^2 s'}$ is a Brownian motion and
$\widetilde{T}_1= \frac{T_r}{r^2}$ is the first time it hits the disk of
radius $1$.
Therefore, from~(\ref{scalingko}), we get the following scaling
relation in distribution for all $r \in\,]0,1]$ under the annealed
measure $\Pb_0$:
%
\begin{eqnarray}\label{invechellsto}
&& (r \eps/2)^{\sfrac{\gamma^2} 2}\int_0^{T_r}
e^{ \gamma\bar{X}_{r
\eps
}(B_s) } \,ds
\nonumber\\[-8pt]\\[-8pt]\nonumber
&&\qquad \mathop{=}^{\mathrm{law}} r^2 e^{\gamma\Omega_r-\sfrac{\gamma^2}{2}
\ln\sklfrac{1}{r}}
(\epsilon/2)^{ \sfrac{\gamma^2} 2}\int_0^{T_1}
e^{
\gamma\bar{X}_{\eps}( B_s)} \,ds.
\end{eqnarray}

From this scaling relation and the above considerations, we deduce that
we can find some variable $N$ with negative moments and such that we
have the following stochastic domination:
\[
Y \geq N(Y_{1}+Y_{2}),
\]
where $(Y_{1},Y_{2})$ are i.i.d. of distribution $Y$, independent of
$N$, $\E[N^{-q}]<+\infty$ (see \cite{AdTay}, Theorem 2.1.1) and $Y$ is
distributed like
$ \lim_{\eps\to0} \epsilon^{\sfrac{\gamma^2} 2} \int_0^{T_1} e^{
\gamma\bar{X}_{\eps}( B_s) } \,ds$. We get~(\ref{equintermediaire})
is true with $r=1$ by adapting \cite{molchan} (see also \cite{Rnew12}, Section~B.4). 
Then one deduces inequality~(\ref{equintermediaire}) from~(\ref
{invechellsto}) for all $r$.\end{pf}

One can then deduce the following results.

\begin{corollary}\label{negative2}
For all $q>0$, there exists some constant $C>0$ (depending on $q$) such
that for all $x\in\R^2$ and all $0\leq s<t \leq1$:
%
\begin{equation}
\label{momnegatiffixe} \sup_{n\geq0} \E_x \biggl[ \biggl(
c_n^{- \sfrac{\gamma^2} 2}\int_s^{t}
e^{ \gamma X_n( B_r) } \,dr \biggr)^{-q} \biggr] \leq C (t-s)^{\xi(-q)}.
\end{equation}
\end{corollary}

\begin{pf} Without loss of generality, we can take $x=0$
and $s=0$ by stationarity of the field $X$ and the strong Markov
property of the Brownian motion.
Then it suffices to show that the expectation in~(\ref{momnegatiffixe})
is finite for $t=1$. Indeed, one can then use the techniques in the
proof of Proposition \ref{lemmamoments} in order to obtain the
right-hand side of~(\ref{momnegatiffixe}) for any $t \in[0,1]$
($x\mapsto x^{-q}$ is convex).

Recall that $T_r$ denote the first exit times of the Brownian motion
out of the ball $B(0,r)$. Recall that, by scaling and \cite{sznit}, Theorem~1.2 (page~93), we have the existence of some absolute constant
$c>0$ such that for all $t>0$
%
\begin{equation}
\label{rappelSzn} \Pb^B_0 (T_{r} \geq t ) \leq
\frac{c t}{r^2} e^{-\vfrac{c t}{r^2}}.
\end{equation}

We have
\[
F(1)^{-q} \leq\sum_{n=1}^{\infty}
\ind_{ \{T_{1/2^n} < 1 \leq
T_{1/2^{n-1}} \}} F({T_{1/2^n})^{-q}}.
\]
%
Therefore, we get by Proposition \ref{negative} and the bound~(\ref
{rappelSzn})
\begin{eqnarray*}
\E_0 \bigl[ F(1)^{-q}\bigr] & \leq&\sum
_{n=1}^{\infty} \E_0\bigl[
\ind_{ \{T_{1/2^n}
< 1 \leq T_{1/2^{n-1}} \}} F({T_{1/2^n})^{-q}}\bigr]
\\
& \leq&\sum_{n=1}^{\infty}
\Pb^B_0 (T_{1/2^{n-1}} \geq1 )^{1/2}
\E_0\bigl[ F(T_{1/2^n})^{-2q} \bigr]^{1/2}
\\
& \leq& C \sum_{n=1}^{\infty} 2^{2n}
e^{-c 2^{2n}} \bigl(2^{n}\bigr)^{2q+\sklvfrac
{q(1+2q)}{2} \gamma^2}.
\end{eqnarray*}
%
This latter series is obviously finite. \end{pf}

\begin{corollary}
Set $\beta=  (1+\frac{\gamma}{2} )^2$. For each $T>0$ and
$\epsilon>0$, there exists a random constant $C>0$ such that $\Pb_x$ a.s.
\[
\forall0\leq s<t\leq T,\qquad\bigl\llvert F(t)-F(s)\bigr\rrvert \geq C\llvert
t-s\rrvert ^{\beta+\epsilon}.
\]
\end{corollary}

\begin{pf} We take $T=1$ for simplicity. Now, we partition
$[0,1]$ into $2^{n}$ dyadic squares $(I_n^j)_{1 \leq j \leq2^{2n}}$ of
equal size. If $p>0$, we get
\begin{eqnarray*}
\Pb_x\biggl( \inf_{1 \leq j \leq2^{n}} F\bigl(I_n^j
\bigr) \leq\frac{1}{2^{(\beta
+\epsilon)n}} \biggr) & \leq&2^{-p(\beta+\epsilon)n} \E_x
\biggl[ \sum_{1 \leq j
\leq2^{n}} F\bigl(I_n^j
\bigr)^{-p} \biggr]
\\
& \leq& C_p 2^{-p(\beta+\epsilon)n} 2^{(1-\xi(-p))n}.
\end{eqnarray*}
We conclude as in the proof of Theorem \ref{app:multform} by choosing
$p>0$ such that $\xi(-p)-1+p(\beta+\epsilon)>0$. \end{pf}

\subsection{Definition of the Liouville Brownian motion}
Before giving the explicit description of the Liouville Brownian
motion, let us first motivate the forthcoming definitions. The reader
may consult \cite{hsu} for an introductory background on Brownian
motion on manifolds. Recall that our purpose is to define the Brownian
motion associated to the metric tensor formally written (using
conventional notation in Riemann geometry) $e^{\gamma X(x)}\,dx^2$, where
$dx^2$ stands for the standard Euclidean metric on $\R^2$. Because of
the obvious divergences of such a direct definition, it is natural to
regularize the field $X$ and to consider instead the Riemann metric tensor
%
\begin{equation}
g_n=c_n^{-\sfrac{\gamma^2}{2}}e^{ \gamma X_n(x)}\,dx^2.
\end{equation}
The renormalization sequence $(c_n^{-\sfrac{\gamma^2}{2}})_n$ appears
here to regulate the divergences when $n\to\infty$.

The Riemann volume on the manifold $(\R^2,g_n)$ is nothing but the
measure $M_n$ defined by~(\ref{def:Mn}). One can also associate to the
Riemann manifold $(\R^2,g_n)$ a Brownian motion $\mathcal{B}^n$.
Following the standard construction, such a Brownian motion can be
constructed as follows. Consider a standard Brownian motion $Z$ on $\R
^2$ and, for any $x\in\R^2$, consider the solution $\LZ^{n,x}$ of the
following stochastic differential equation:
%
\begin{equation}
\label{def.LBM1}
\cases{
\LZ^{n,x}_{\mathbf{t}=0}=x,
\vspace*{3pt}\cr
\displaystyle d\LZ^{n,x}_\mathbf{t}= c_n^{\sfrac{\gamma^2}{4}}e^{-\sfrac{\gamma}
2 X_n(\LZ
^{n,x}_\mathbf{t})}
\,dZ_\mathbf{t}.}
\end{equation}
Notice that it is not clear for the time being that the solution of
this SDE has an infinite lifetime.
Actually, we will see that the solution of such a SDE does not converge
in probability a $n\to\infty$. Yet, it is enough to investigate
convergence in law in the space of continuous functions $C(\R_+,\R^2)$.
By using the Dambis--Schwarz theorem, the solution of~(\ref{def.LBM1})
as the same law as the following process.

\begin{definition}\label{d.elbm2} For any $n\geq1$, we define
%
\begin{equation}
\LB^{n}_\mathbf{t}= B_{\langle\LB^{n} \rangle_\mathbf{t}},
\end{equation}
where $(B_r)_{r\geq0}$ is a planar Brownian motion, independent of the
MFF $X$ and where the quadratic variation $\langle\LB^{n} \rangle$
of $\LB^{n,x}$
is defined as follows:
%
\begin{equation}
\label{eq.qvdambis} \bigl\langle\LB^{n} \bigr\rangle_\mathbf{t}:=
\inf\biggl\{s\geq0: c_n^{-\sfrac{\gamma^2}{2}}\int_{0}^s
e^{ \gamma X_n(B_u)} \,du \geq\mathbf{t}\biggr\}.
\end{equation}
We call this process $n$-\textit{regularized Liouville Brownian motion}
($n$-LBM for short).
\end{definition}

The infimum in~(\ref{eq.qvdambis}) is necessarily finite. Indeed we have
\[
c_n^{-\sfrac{\gamma^2}{2}}\int_{0}^s
e^{ \gamma X_n(B_u)} \,du\geq c_n^{-\sfrac{\gamma^2}{2}}e^{\min_{x\in B(0,1)}X_n(x)} \int
_0^t\ind _{\{
B_r\in B(0,1)\}} \,dr
\]
and the latter quantity tends to $\infty$ when $t\to\infty$ via
recurrence arguments.


Observe that~(\ref{eq.qvdambis}) amounts to saying that the increasing
process $\langle\mathcal{B}^{n,x}\rangle:\R_+\to\R_+$ is just the
inverse function of the PCAF
%
\begin{equation}
\label{fnt} F_n(t)=c_n^{-\sfrac{\gamma^2}{2}}\int
_0^te^{\gamma X_n(B_r)} \,dr
\end{equation}
studied in Section~\ref{sub:pcaf}. In particular, we can assert now
that the solution of~(\ref{def.LBM1}) has infinite lifetime because
$F_n(t)\to\infty$ as $t\to\infty$ (adapt, e.g., the argument of
item 4 in Theorem \ref{th:PCAF}). Several standard facts can be deduced
from the smoothness of $X_n$.

%
\begin{proposition}\label{prop.fellern}
Let $n\geq1$ be fixed. $\Pb^X$-a.s., the $n$-regularized Liouville
Brownian motion $\mathcal{B}^n$ induces a Feller diffusion on $\R^2$.
Let us denote by $(P_\mathbf{t}^{\gamma,n})_{\mathbf{t}\geq0}$ its
semigroup. Also,
$\mathcal{B}^n$ is reversible with respect to the Riemann volume $M_n$,
which is therefore invariant for $\mathcal{B}^n$.
\end{proposition}

For a proof of this proposition, the reader may consult \cite{fuku}, for
instance, though the Feller property is not discussed but this is easy
(for the same reason as the Liouville Brownian motion is Feller; see
item 1 of Proposition \ref{referee} below).

As we have seen that the PCAFs $(F_n)_n$ converge toward $F$, it is
natural to introduce the following.

\begin{definition}\label{def:LBM}
Assume $\gamma<2$. $\Pb^X$ almost surely, we define the \textit{Liouville
Brownian motion} (LBM for short) as the time changed Brownian motion
\[
\mathcal{B}_t=B_{F^{-1}(t)},
\]
where $F$ is the PCAF of Theorem \ref{th:PCAF}.
\end{definition}

Observe that the continuity of $F$ makes sure that the LBM does not
get stuck in some area of the state space $\R^2$. Typically, this
situation may happen over areas where the field $X$ takes large values,
therefore, having as consequence to slow down the LBM. Furthermore,
strict monotonicity of $F$ ensures that the LBM possesses no jumps and
the fact that $F$ tends to $\infty$ as $\mathbf{t}\to\infty$ makes
sure that
the LBM has an infinite lifetime, it does not ``reach $\infty$'' in
finite time.

\subsection{Main properties of the LBM}\label{ss.prop}
In this subsection, we collect a few important properties of the LBM
that can be deduced from our previous analysis. The reader is referred
to \cite{fuku} for the classical terminology on Dirichlet forms.

\begin{theorem}\label{th:sgp}
For $\gamma\in[0,2[$, $\Pb^X$-a.s., the LBM is a strong Markov process
with continuous sample paths, that is a diffusion process. If we denote
by $(P^{\gamma,X}_\mathbf{t})_{\mathbf{t}\geq0}$ the associated
semigroup on $L^2(\R
^2,M)$ then:
\begin{longlist}[2.]
\item[1.]$(P^{\gamma,X}_\mathbf{t})_{\mathbf{t}\geq0}$ is strongly
continuous on $L^2(\R
^2,M)$ and symmetric.
\item[2.] the associated Dirichlet form is strongly local, regular and
possesses $C^\infty_c(\R^2)$ as special core.
\end{longlist}
\end{theorem}

Observe that the semigroup $(P^{\gamma,X}_\mathbf{t})_{\mathbf
{t}\geq0}$ is itself
random as it depends on the randomness of the MFF $X$. Furthermore, the
above theorem entails that we can classically extend the semigroup
$(P^{\gamma,X}_\mathbf{t})_{\mathbf{t}\geq0}$ to a
strongly-continuous semi-group on
$L^p(\R^2, M)$ for all $1\leq p <\infty$.

\begin{pf*}{Proof of Theorem \ref{th:sgp}}
From \cite{fuku}, Theorem~6.1.1, the LBM is a Hunt process, therefore strong Markov. Item~1 and 2 then result from \cite{fuku}, Theorem~6.2.1. Continuity of
sample paths results from the fact that $F$ is strictly increasing
(Theorem \ref{th:PCAF} item 3), in which case $F^{-1}$ is
continuous.\end{pf*}

\begin{proposition}\label{referee}
For $\gamma\in[0,2[$, $\Pb^X$-a.s.:
\begin{longlist}[2.]
\item[1.] the law in $C(\R_+)$ of the LBM $\LB$ under $\Pb^B_x$ is
continuous with respect to $x$.
\item[2.] for all $x\in\R^2$, under $\Pb^B_x$ the $n$-LBM $(\LB^n)_n$
converges in law in $C(\R_+)$ equipped with the topology of uniform
convergence over compact sets toward $\LB$,\vspace*{1pt}
\item[3.] the semigroup $(P^{\gamma,X}_\mathbf{t})_{\mathbf{t}\geq0}$
is the limit of the
semigroups $(P^{\gamma,n})_n$ as $n\to\infty$ in the sense that
\[
\forall f\in C_b\bigl(\R^2\bigr),\qquad\lim
_{n\to\infty}P^{\gamma,n}_\mathbf{t} f(x)=P^{\gamma,X}_\mathbf{t}f(x),
\]
\item[4.] for all $x\in\R^2$ and for all $z\in\R^2$, $\Pb^B_x
(\liminf_{t\to\infty} \llvert  \LB_t - z\rrvert  =0 )=1$,
\item[5.] for all $x\in\R^2$, $\Pb^B_x (\limsup_{t\to\infty} \llvert  \LB_t\rrvert
=\infty )=1$.
\end{longlist}
\end{proposition}

\begin{pf} Let us prove item 1. Observe that the mapping
\[
(w,v)\in C\bigl(\R_+,\R^2\bigr)\times C(\R_+,\R_+)\mapsto w\circ
v^{-1},
\]
where
$v^{-1}(t)=\inf\{s\geq0; v(s)>t\}$ is continuous at all these pairs
$(w,v)$ such that $v$ is strictly increasing and $\lim_{s\to\infty
}v(s)=+\infty$. Therefore, item 1 results from Theorem \ref{th:PCAF}
items 3${}+{}$4${}+{}$5.

The same argument and Theorem \ref{th:PCAF} items $2+3+4$ prove item 2.
Item 3 is a direct consequence of item 2. Items 4 and 5 results from the
fact a standard Brownian motion satisfies items $4+5$ and the fact that
$\LB$ is a time changed Brownian motion, with a continuous time change
that goes to $\infty$ as $t\to\infty$. \end{pf}

\begin{corollary}
For $\gamma\in[0,2[$, $\Pb^X$-a.s., the semigroup $(P^{\gamma
,X}_\mathbf{t}
)_{\mathbf{t}\geq0}$ is Feller.
\end{corollary}

For $\gamma<2$, we define the \textit{Liouville Laplacian} $\Delta_X$ as
the generator of the Liouville Brownian motion times the usual extra
factor $\sqrt{2}$. The Liouville Laplacian corresponds to an operator
which can formally be written as
\[
\Delta_X=e^{-\gamma X (x)}\Delta
\]
and can be thought of as the Laplace--Beltrami operator of
$2d$-Liouville quantum gravity.

\subsection{Asymptotic independence of the Liouville Brownian motion
and the Euclidean Brownian motion}
Recall the SDE~(\ref{def.LBM1}) with solution $\LZ^{n,x}$, which is
measurable with respect to the planar Brownian motion $Z$. In this
section, we prove that $(\LZ^{n,x})_n$ does not converge in probability
as $n\to\infty$. This will show that the time change representation of
the LBM is more relevant.

\begin{theorem}\label{th.asympind}
If $\gamma<2$, $\Pb^X$-a.s. and for all $x\in\R^2$, the pair of
processes $(Z,\LZ^{n,x})_n$ converges in law toward a couple $(Z,\LZ
^x)$. The planar Brownian motion $Z$ and the LBM $\LZ^x$ are independent.
\end{theorem}

The above theorem shows that some extra randomness is created by taking
the limit $n\to\infty$. Indeed, the $n$-LBM is a measurable function of
the planar Brownian motion $Z$. Yet, Liouville/Euclidean Brownian
motions are independent at the limit, showing that convergence in
probability cannot hold.

\begin{pf*}{Proof of Theorem \ref{th.asympind}} Before beginning
the proof, let us first clarify a few points. The $n$-LBM~(\ref
{def.LBM1}) involves the planar Brownian motion $Z$. An equivalent
definition in law of this $n$-LBM is given in Definition \ref{d.elbm2}
by means of another Brownian motion $B$, constructed via the
Dambis--Schwarz theorem. As such, it implicitly depends on $n$ as well
as $Z$. It is therefore relevant to write explicitly this dependence in
this proof. So we will write $B^n$ instead of $B$.

We begin with writing explicitly the dependence between $Z$ and $B^n$.
The Dambis--Schwarz theorem tells us that [recall~(\ref{fnt})]
\[
B^{n}_t=c_n^{\sfrac{\gamma^2} 4 }\int
_0^{F_n(t)} e^{-\sfrac{\gamma}
2 X_n(\LZ^{n,x}_\u) } \,dZ_\u=
\LZ^{n,x}_{F_n(t,x)}.
\]

Now we prove the asymptotic independence of $Z$ and $B^n$. Let us
compute their predictable bracket by making the change of variables:
\begin{eqnarray*}
\bigl\langle B^{n},Z\bigr\rangle_t&=&c_n^{\sfrac{\gamma^2} 4 }
\int_0^te^{-\sklfrac {\gamma} 2 X_n(\LZ^{n,x}_r)} \,dr
\\
& =&c_n^{\sfrac{\gamma^2}{4} }\int_0^t
e^{-\sklfrac{\gamma
}{2}X_n(B^n_{F_n^{-1}(r) } )} \,dr
\\
& =&c_n^{-\sfrac{\gamma^2}{4} }\int_0^{ F_n(t)}e^{ \sklfrac{\gamma} 2
X_n(B^n_{u})}
\,du
\\
&=&c_n^{-\sfrac{\gamma^2} 8 }\times c_n^{-\sfrac{\gamma^2} 8 }\int
_0^{F_n(t)}e^{ \sklfrac{\gamma} 2 X_n(B_{u})} \,du.
\end{eqnarray*}
From Theorem \ref{th:PCAF}, the above expression corresponds to
multiplying a tight family by a factor $c_n^{-\sfrac{\gamma^2} 8 }$
that goes to $0$ as $n\to\infty$. Therefore, $\Pb^X$-a.s., the family
$(\langle B^{n},Z\rangle)_n$ (as random functions of $t$) converges
$\Pb
^Z$-a.s. in $C(\R_+)$ toward $0$. The pair $(B^{n},Z)$ therefore
converges in law in $C(\R_+)$ toward a pair $(B,Z)$ of Brownian
motions, the brackets of which vanish. Knight's theorem \cite{karatzas}, Theorem~4.13, implies that $B$ and $Z$ are independent (see also the
Appendix in \cite{pitman}). As a measurable function of $B$, the
Liouville Brownian motion is independent of $Z$.\end{pf*}

\subsection{Liouville Brownian motion defined on other geometries:
Torus, sphere and planar domains}

So far, we constructed in detail the LBM for the (Massive) Free Field
on $\R^2$. Our method applies to other two-dimensional manifolds like
the torus, sphere, planar domains$,\ldots$ equipped with a log-correlated
Gaussian field (of special interest is the case of Gaussian Free
Field), stationary or not. The main reason is that Kahane's theory
remains valid on $C^1$-manifolds (see \cite{cfKah,review}).
Intuitively, this is just because such manifolds are locally isometric
to open domains of the Euclidean space.

There is at least one point in our proofs that must be changed in order
to apply to the torus or the sphere, or any compact manifold without
boundary: the fact that $\lim_{t\to\infty}F(x,t)=+\infty$. Indeed, our
proof uses the ``infinite volume'' of the plane. In the case of the torus
or sphere, the strategy is much simpler because of compactness
arguments: the standard Brownian motion on $\S^2$ or $\T^2$ possesses
an invariant probability measure, call it $\mu$, which is nothing but
the volume form of $\S^2$ or $\T^2$. Apply the ergodic theorem to prove
that $\Pb^X\otimes\Pb^B_\mu$ almost surely:
\[
\lim_{t\to\infty}\frac{F(t)}{t}=G,
\]
for some random variable $G$, which is shift-invariant. Since the
Brownian motion on the sphere is ergodic, $G$ is measurable with
respect to the sigma algebra generated by $\sigma\{X_x;x\in\T^2
\mbox{ or }\S^2\}$. It is not clear that $G$ is constant. Yet, the set
$\{G>0\}$ is measurable with respect to the asymptotic sigma-algebra of
the $(Y_n)_n$. Therefore, $\Pb^\mu$ almost surely, the set $\{G>0\}$
has $\Pb^X$-probability $0$ or $1$. Since $G$ has expectation $1$, this
set has $\Pb^X$-probability $1$.
Therefore, $\Pb^X$ almost surely, the change of times $F$ goes to
$\infty$ as $t\to\infty$ for $\mu$ almost every $x$. Then use the
coupling trick to deduce that the property holds for all starting points.

\begin{appendix}
\section{Kahane's convexity inequality}\label{sec.chaos}
For the classical terminology of Gaussian multiplicative chaos, the
reader is referred to \cite{cfKah} (see also \cite{review}).

We consider a locally compact separable metric space $(D,d)$, a Radon
measure $\nu$ on the Borel subsets of $(D,d)$ and two Gaussian random
distributions $Y,Y'$ (in the sense of Schwartz) with respective
covariance kernels $K,K'$, which are of $\sigma$-positive type. We
assume that the Gaussian multiplicative chaos associated to $(Y,\nu)$
and $(Y',\nu)$ are strongly nondegenerate (e.g., the kernels of
$Y_n$ [see~(\ref{kernelY})] or the kernel $\gamma^2 G_m$ of~(\ref
{MFF1}) for $0\leq\gamma^2<4$).
Here, we recall the following standard lemma that can be found in \cite
{cfKah}.

\begin{lemma}\label{cvx}
Let $F:\R_+\to\R$ be some convex function such that
\[
\forall x\in\R_+,\qquad\bigl\llvert F(x)\bigr\rrvert \leq M\bigl(1+\llvert x
\rrvert ^\beta\bigr),
\]
for some positive constants $M,\beta$.

\begin{longlist}[(2)]
\item[(1)] Assume that $ K(u,v)\leq K'(u,v)$ for all $u,v\in D$. Then
\[
\E \biggl[F \biggl(\int_De^{Y_r-\sklfrac{1}{2}\E[Y_r^2] } \nu(dr)
\biggr) \biggr]\leq \E \biggl[F \biggl(\int_De^{Y'_r-\sklfrac{1}{2}\E[Y_r^{\prime 2}] }
\nu (dr) \biggr) \biggr].
\]

\item[(2)] If $K(u,v)\leq K'(u,v)+ C$ for some constant $C>0$ for all $u,v\in
D$ then
\[
\E \biggl[F \biggl(\int_De^{Y_r-\sklfrac{1}{2}\E[Y_r^2] } \nu(dr)
\biggr) \biggr]\leq\E \biggl[F \biggl(e^{\sqrt{C}Z-C/2}\int_De^{Y'_r-\sklfrac{1}{2}\E
[Y_r^{\prime 2}] }
\nu(dr) \biggr) \biggr],
\]
where $Z$ is a standard Gaussian random variable
independent of the other random quantities.
\end{longlist}
\end{lemma}

\section{Finiteness of the moments}\label{app:moments}
In this section, our only goal is to prove that
\[
\forall x\in\R^2,\qquad\E_x\bigl[F(t)^p
\bigr]<+\infty
\]
for $p\in[0,4/\gamma^2[$. We only treat the case when $2\leq\gamma
^2<4$ and, therefore, $1<p<4/\gamma^2$ (hence $p<2$). This is
mathematically the most complicated part and notationally the easiest
part. The case $0\leq\gamma^2<2$ is discussed in Remark \ref{rem:mom}.

Furthermore, by stationarity of the field $X$, we may assume that
$x=0$. By using the concavity of the mapping $x\mapsto x^{p/2}$ and the
Jensen inequality, we get
%
\begin{eqnarray}\label{mainmaj}
\E_0\bigl[F(t)^p \bigr] &\leq&
\E^X \bigl[\E^B_0\bigl[ F(t)^2
\bigr]^{p/2} \bigr] 
\nonumber\\[-8pt]\\[-8pt]\nonumber
&\leq&
\E^X \biggl[ \biggl(\int_{\R^2}\int
_{\R^2}f(x,y)M(dx)M(dy) \biggr)^{p/2} \biggr],
\end{eqnarray}
where we have set
%
\begin{equation}
\label{def:f} f(x,y)= \int_0^t\int
_s^t e^{-\afrac{\llvert  x\rrvert  ^2}{2s}-\afrac
{\llvert  y-x \rrvert  ^2}{2\llvert  r-s\rrvert  }} \frac{dr\,ds}{4\pi^2s\llvert  r-s\rrvert  }.
\end{equation}
So we just have to prove that the expectation in the above right-hand
side is finite.
In what follows, $\xi_M$ stands for the structure exponent of the
measure $M$. Recall that, in dimension $d$, it reads
%
\begin{equation}
\label{def:xid} \xi_M(p)=\biggl(d+\frac{\gamma^2}{2}\biggr)p-
\frac{\gamma^2}{2}p^2.
\end{equation}
Of course, we can take here $d=2$. But is worth recalling this fact
since it will happen that some arguments below will be carried out in
dimension $1$. So the reader will take care of replacing $d$ by $1$
when reading a proof in dimension $1$. The main idea of our proof is
the following. First, we observe that the function $f$ possesses
singularities. They are logarithmic (see below) when $x$ or $\llvert  x-y\rrvert  $ is
close to $0$. We will also have to treat the behavior near infinity. So
we split the space $\R^2\times\R^2$ into $3$ domains
%
\begin{eqnarray}
\mathcal{D}_1&=&\bigl\{\llvert x\rrvert \leq1, \llvert x-y\rrvert
\leq1\bigr\},\qquad\mathcal{D}_2=\bigl\{ \llvert x\rrvert \geq 1,
\llvert x-y\rrvert \leq1\bigr\},
\nonumber\\[-8pt]\\[-8pt]\nonumber
\mathcal{D}_3&=&\bigl\{ \llvert
x-y\rrvert \geq1\bigr\}
\end{eqnarray}
and, by subadditivity of the mapping $x\in\R_+\mapsto x^{p/2}$, it
suffices to evaluate the quantity in the right-hand side of~(\ref
{mainmaj}) on each of these three domains.

Concerning the behavior of $f$, we claim the following.

%
\begin{lemma}\label{lem:f}
We have:
\begin{longlist}[2.]
\item[1.] for all $x,y\in\R^2$: $f(x,y)\leq D(1+\ln_+\frac
{1}{\llvert  x-y\rrvert  })(1+\ln
_+\frac{1}{\llvert  x\rrvert  })$,
\item[2.] for all $\llvert  x\rrvert  \geq1$ and $\llvert  x-y\rrvert  \leq1$: $f(x,y)\leq D(1+\ln
_+\frac
{1}{\llvert  x-y\rrvert  })\exp (-\frac{\llvert  x\rrvert  ^2}{4t} )$, for some constant
$D>0$.
\end{longlist}
\end{lemma}

\begin{pf} Recall~(\ref{def:f}). By making successive
changes of variables, we obtain
\begin{eqnarray*}
f(x,y)&=&\int_0^t\int_0^{\vfrac{t-s}{\llvert  x-y\rrvert  ^2}}
e^{-\afrac
{\llvert  x\rrvert  ^2}{2s}-\afrac
{1}{2r}} \frac{dr\,ds}{4\pi^2sr}
\\
&=& \int_0^{\sfrac{t}{\llvert  x\rrvert  ^2}}
\int_0^{\vfrac
{t-s\llvert  x\rrvert  ^2}{\llvert  x-y\rrvert  ^2}} e^{-\afrac{1}{2s}-\afrac{1}{2r}} \frac
{dr\,ds}{4\pi
^2sr}
\\
&\leq& g\biggl(\frac{t}{\llvert  x\rrvert  ^2}\biggr)g\biggl(\frac{t}{\llvert  x-y\rrvert  ^2}\biggr),
\end{eqnarray*}
where we have set
\[
g(t)=\int_0^te^{-\afrac{1}{2s}}
\frac{ds}{2\pi s }.
\]
It is obvious to check that, for some constant $D>0$, we have $g(t)\leq
D(1+\ln_+ t)$, which completes the proof of item 1. The proof of item 2
is similar and left to the reader.\end{pf}

\begin{notation}
Until the end of the proof, we will only deal with expectations with
respect to the measure $M$. So there is no need to keep on using the
superscript $X$ of $\E^X$ and we will just write $\E$ instead of $\E^X$.
\end{notation}

\subsection*{Domain $\{\llvert x\rrvert\leq1,\llvert x-y\rrvert\leq1\}$}
The main purpose of this part is to show that
%
\begin{equation}
\label{lem:Mplus} \E \biggl[ \biggl(\int_{\max(\llvert  x\rrvert,\llvert  y\rrvert  )\leq1}
\frac{1}{\llvert  x\rrvert  ^\delta
\llvert  x-y\rrvert  ^\delta
}M(dx)M(dy) \biggr)^{p/2} \biggr]<+\infty.
\end{equation}
%

The first step is to prove the following.

\begin{lemma}\label{lem:Mzero}
For $\gamma^2<4$ and $p\in\,]1,\frac{4}{\gamma^2}[$, there exist
$\delta
>0$ and $C>0$ such that
for all $n\geq0$,
\[
\E \biggl[ \biggl(\int_{ \max(\llvert  x\rrvert,\llvert  y\rrvert  )\leq2^{-n}} \frac{1}{\llvert  x\rrvert  ^\delta
}M(dx)M(dy)
\biggr)^{p/2} \biggr]\leq C2^{-n(\xi_M(p)-\vfrac{\delta p}{2})}.
\]
\end{lemma}

\begin{pf} We carry out the proof in dimension $1$ since,
apart from notational issues, the dimension $2$ does not raise any
further difficulty. We first have to prove
\[
\E \biggl[ \biggl(\int_{(x,y)\in[0,1]^2} \frac{1}{\llvert  x\rrvert  ^\delta} M(dx)M(dy)
\biggr)^{p/2} \biggr]<+\infty.
\]
Furthermore, from Kahane's convexity inequalities, it suffices to prove
the above lemma for any $1d$ log-correlated Gaussian field. Let us use
the kernel\vadjust{\goodbreak} $K(x,y)=\ln_+\frac{1}{\llvert  x-y\rrvert  }$ of \cite{bacry}. In fact,
whatever the covariance kernel, we just need to use the property~(\ref
{eq:pls}), which is shared by all the reasonable $1d$ log-correlated
Gaussian fields (see \cite{review} for more on this).

We also remind the reader that the above integral is finite for $\delta
=0$ (see \cite{cfKah}). Therefore, by using subadditivity of the
mapping $x\mapsto x^{p/2}$, we have
\begin{eqnarray*}
&& \E \biggl[ \biggl(\int_{ [0,1]^2} \frac{1}{\llvert  x\rrvert  ^\delta} M(dx)M(dy)
\biggr)^{p/2} \biggr]
\\
&&\qquad =\E \Biggl[ \Biggl(\sum
_{n=0}^{\infty}\int_{ [2^{-n-1}, 2^{-n}]\times[0,1]}
\frac{1}{\llvert  x\rrvert  ^\delta} M(dx)M(dy) \Biggr)^{p/2} \Biggr]
\\
&&\qquad \leq\sum_{n=0}^{\infty}\E \biggl[ \biggl(\int
_{[2^{-n-1},
2^{-n}]\times
[0,1]} \frac{1}{\llvert  x\rrvert  ^\delta} M(dx)M(dy) \biggr)^{p/2}
\biggr]
\\
&&\qquad \leq\sum_{n=0}^{\infty}2^{\delta(n+1)p/2}\E
\bigl[ \bigl( M\bigl(\bigl[2^{-n-1},2^{-n}\bigr]\bigr)M
\bigl([0,1]\bigr) \bigr)^{p/2} \bigr].
\end{eqnarray*}
Now we use the standard inequality $ab\leq\epsilon a^2+\frac
{b^2}{\epsilon}$ for any $\epsilon>0$ and subadditivity of the mapping
$x\mapsto x^{p/2}$ to get $(ab)^{p/2}\leq\epsilon^{p/2}a^p+\epsilon
^{-p/2}b^p$. Therefore, with $a=M([2^{-n-1},2^{-n}])$, $b=M([0,1])$ and
$\epsilon=2^{(n+1)\xi_M(p)/p}$, we obtain
\begin{eqnarray*}
\E \bigl[ \bigl( M\bigl(\bigl[2^{-n-1},2^{-n}\bigr]\bigr)M
\bigl([0,1]\bigr) \bigr)^{p/2} \bigr]&\leq &\bigl(2^{(n+1)\xi_M(p)/p}
\bigr)^{p/2}\E \bigl[ M\bigl(\bigl[2^{-n-1},2^{-n}
\bigr]\bigr)^{p} \bigr]
\\
&&{}+\bigl(2^{(n+1)\xi_M(p)/p}\bigr)^{-p/2}\E \bigl[ M\bigl([0,1]
\bigr)^{p} \bigr].
\end{eqnarray*}
By using~(\ref{eq:pls}), we get
\[
\E \bigl[ M\bigl(\bigl[2^{-n-1},2^{-n}\bigr]
\bigr)^{p} \bigr]\leq C_p 2^{-(n+1)\xi_M(p)}
\]
and plugging this relation into the above expression yields:
\begin{eqnarray*}
\E \bigl[ \bigl( M\bigl(\bigl[2^{-n-1},2^{-n}\bigr]\bigr)M
\bigl([0,1]\bigr) \bigr)^{p/2} \bigr]&\leq& 2^{-(n+1)\xi_M(p)/2}
\bigl(C_p+\E \bigl[ M\bigl([0,1]\bigr)^{p} \bigr] \bigr).
\end{eqnarray*}
To sum up, we have
\begin{eqnarray*}
&&\E \biggl[ \biggl(\int_{(x,y)\in[0,1]^2} \frac{1}{\llvert  x\rrvert  ^\delta} M(dx)M(dy)
\biggr)^{p/2} \biggr]
\\
&&\qquad  \leq \bigl(C_p+\E \bigl[ M\bigl([0,1]
\bigr)^{p} \bigr] \bigr)\sum_{n=0}^{\infty}2^{-(n+1)(\xi_M(p)/2-\delta p/2)}.
\end{eqnarray*}
So, $\delta$ can clearly be chosen small enough to make the above
series convergent.

Once the finiteness of the expectation is proved, the statement results
from a scaling argument. For $\lambda<1$, the measure $M$ satisfies
(see \cite{bacry} but this is elementary) the following relation in law:
%
\begin{equation}
\label{esi} \bigl(M(\lambda A)\bigr)_{A\subset[0,1]} =\bigl(
\lambda^{1+\sfrac{\gamma^2}{2}} e^{\gamma\Omega_\lambda}M(A)\bigr)_{A\subset[0,1]},
\end{equation}
where $\Omega_\lambda$ is a centered Gaussian random variable with
variance $-\ln\lambda$ independent of $(M(A))_{A\subset[0,1]}$. Thus,
we have
\begin{eqnarray*}
&& \E \biggl[ \biggl(\int_{ [0,\lambda]^2} \frac{1}{\llvert  x\rrvert  ^\delta} M(dx)M(dy)
\biggr)^{p/2} \biggr]
\\
&&\qquad  =\lambda^{p(1+\sfrac{\gamma
^2}{2})-\delta p/2}\E\bigl[e^{p \gamma\Omega_\lambda}\bigr] \E
\biggl[ \biggl(\int_{
[0,1]^2} \frac{1}{\llvert  x\rrvert  ^\delta} M(dx)M(dy)
\biggr)^{p/2} \biggr].
\end{eqnarray*}
The result follows by taking $\lambda=2^{-n}$.\end{pf}

\begin{lemma}\label{lem:Mdiag}
For any $\gamma^2<4$ and $p\in\,]1,\frac{4}{\gamma^2}[$, there exist
$\delta>0$ and a constant $C>0$ 
such that for all $n$:
\[
\E \biggl[ \biggl(\mathop{\int_{\max(\llvert  x\rrvert,\llvert  y\rrvert  )\leq1}}_{{2^{-n-1}\leq\llvert  x-y\rrvert  \leq
2^{-n}}}
\frac{1}{\llvert  x\rrvert  ^\delta}M(dx)M(dy) \biggr)^{p/2} \biggr]\leq
\frac
{C}{1-\delta p/2}2^{-n(\xi_M(p)-2)}.
\]
\end{lemma}

\begin{pf} Once again and for the same reason as
previously, we carry out the proof in dimension $1$. In that case, we
have to prove
\[
\E \biggl[ \biggl(\mathop{\int_{(x,y)\in[0,1]^2}}_{{2^{-n-1}\leq
\llvert  x-y\rrvert  \leq
2^{-n}}}
\frac{1}{\llvert  x\rrvert  ^\delta} M(dx)M(dy) \biggr)^{p/2} \biggr]\leq
C2^{-n(\xi
_M(p)-1)}.
\]
Once again Kahane's convexity inequality shows that we can take the
kernel $K(x,y)=\ln_+\frac{1}{\llvert  x-y\rrvert  }$ of \cite{bacry}. We will use the
following elementary geometric argument: for any $n\geq1$, the set of
points $2^{-n}$-close to the diagonal
\[
\bigl\{(x,y)\in[0,1]^2;\llvert x-y\rrvert \leq2^{-n}\bigr
\}
\]
is entirely recovered by the union for $k=0,\dots,2^n-2$ of the
(overlapping) squares $[\frac{k}{2^n},\frac{k+2}{2^n}]^2$.
Therefore, by using subadditivity of the mapping $x\mapsto x^{p/2}$,
we have
\begin{eqnarray*}
&& \E \biggl[ \biggl( \mathop{\int_{x,y\in[0,1]}}_{{2^{-n-1}\leq\llvert  x-y\rrvert  \leq
2^{-n}}}
\frac{1}{\llvert  x\rrvert  ^\delta} M(dx)M(dy) \biggr)^{ \sfrac{p}{2}} \biggr] %
\\
&&\qquad  \leq\sum_{k=0,\dots,2^n-2} \E \biggl[ \biggl( \int
_{ x,y\in[\sfrac
{k}{2^n},\vfrac{k+2}{2^n}]^2} \frac{1}{\llvert  x\rrvert  ^\delta} M(dx)M(dy) \biggr)^{\sfrac
{p}{2}}
\biggr]
\\
&&\qquad  \leq\E \biggl[ \biggl( \int_{ x,y\in[0,2^{-n+1}]^2} \frac
{1}{\llvert  x\rrvert  ^\delta
}
M(dx)M(dy) \biggr)^{\sfrac{p}{2}} \biggr]
\\
&&\quad\qquad{} +\sum_{k=1,\dots,2^n-2}\frac{2^{n \delta p/2}}{k^{\delta p/2}} \E \biggl[ M
\biggl(\biggl[\frac{k}{2^n},\frac{k+2}{2^n}\biggr] \biggr)^{p}
\biggr].
\end{eqnarray*}
By stationarity and scale invariance~(\ref{esi}), we get
\begin{eqnarray*}
&& \sum_{k=1,\dots,2^n-2}\frac{2^{n \delta p/2}}{k^{\delta p/2}} \E \biggl[ M \biggl(
\biggl[\frac{k}{2^n},\frac{k+2}{2^n}\biggr] \biggr)^{p}
\biggr]
\\
&&\qquad  \leq 2^{n \delta p/2}\sum_{k=1,\dots,2^n-2}
\frac
{1}{k^{\delta p/2}}\E \bigl[ M \bigl(\bigl[0,2^{-n+1}\bigr]
\bigr)^{p} \bigr]
\\
&&\qquad  \leq 2^{n \delta p/2}2^{-(n-1)\xi_M(p)}\E \bigl[ M \bigl([0,1]
\bigr)^{p} \bigr]\sum_{k=1,\dots,2^n-2}
\frac{1}{k^{\delta
p/2}}
\\
&&\qquad  \leq\frac{C}{1-\delta p/2}2^{-n(\xi_M(p)-1)},
\end{eqnarray*}
where $C$ only depends on $\E [ M ([0,1] )^{p} ]$. We
conclude with Lemma \ref{lem:Mzero} provided we impose $\delta
p/2<1$.\end{pf}

Now we are equipped to prove~(\ref{lem:Mplus}). Choose another $\delta
'>0$ such that $0<\delta' <\frac{2(\xi_M(p)-2)}{p}$. By using Lemma
\ref
{lem:Mdiag} and by subadditivity, we have
%
\begin{eqnarray}\label{om}
&& \E \biggl[ \biggl(\int_{\max(\llvert  x\rrvert,\llvert  y\rrvert  )\leq1}\frac{1}{\llvert  x\rrvert  ^\delta
\llvert  x-y\rrvert  ^{\delta
'}}M(dx)M(dy)
\biggr)^{p/2} \biggr]
\nonumber
\\
&&\qquad \leq\sum_{n=0}^{+\infty}\E \biggl[ \biggl(
\mathop{\int_{\max
(\llvert  x\rrvert,\llvert  y\rrvert  )\leq
1}}_{{2^{-n-1}\leq\llvert  x-y\rrvert  \leq2^{-n}}}\frac{1}{\llvert  x\rrvert  ^\delta\llvert  x-y\rrvert  ^{\delta
'}}M(dx)M(dy)
\biggr)^{p/2} \biggr]
\nonumber\\[-8pt]\\[-8pt]\nonumber
&&\qquad \leq\sum_{n=0}^{+\infty}2^{(n+1)\vfrac{\delta' p}{2}}\E
\biggl[ \biggl(\mathop{\int_{\max(\llvert  x\rrvert,\llvert  y\rrvert  )\leq1}}_{{2^{-n-1}\leq\llvert  x-y\rrvert  \leq2^{-n}}}
\frac
{1}{\llvert  x\rrvert  ^\delta}M(dx)M(dy) \biggr)^{p/2} \biggr]
\nonumber
\\
&&\qquad \leq\sum_{n=0}^{+\infty}2^{(n+1)\vfrac{\delta' p}{2}}C2^{-n(\xi
_M(p)-2)}.\nonumber
\end{eqnarray}
Since the latter series converges, the proof of~(\ref{lem:Mplus}) is complete.

By gathering~(\ref{lem:Mplus}) and Lemma \ref{lem:f} item 1, we deduce
%
\begin{equation}
\label{cor:M1} \E \biggl[ \biggl(\int_{\llvert  x\rrvert  \leq1,\llvert  x-y\rrvert  \leq1}f(x,y)M(dx)M(dy)
\biggr)^{p/2} \biggr]<+\infty.
\end{equation}

\subsection*{Domain $\{\llvert x-y\rrvert\geq1\}$}
Let us now investigate the situation when $\llvert  x-y\rrvert  \geq1$. This is the
easy part because, in that case, the measures $M(dx)$ and $M(dy)$ are
``almost'' independent. Therefore, we can proceed more directly in the
computations.\vadjust{\goodbreak} We use the Jensen inequality with the concave function
$x\mapsto x^{p/2}$ to get
\begin{eqnarray*}
&& \E \biggl[ \biggl( \int_{\llvert  x-y\rrvert  \geq1}f(x,y)M(dx)M(dy)
\biggr)^{p/2} \biggr]
\\
&&\qquad \leq \biggl(\E \biggl[\int_{\llvert  x-y\rrvert  \geq1}f(x,y)M(dx)M(dy)
\biggr] \biggr)^{p/2}
\\
&&\qquad \leq \biggl( \int_{\llvert  x-y\rrvert  \geq1}f(x,y)e^{\gamma^2G_m(x,y)} \,dx\,dy
\biggr)^{p/2}.
\end{eqnarray*}
Since $\llvert  x-y\rrvert  \geq1$, we have $G_m(x,y)\leq C$ for some fixed positive
constant $C$. We deduce that the above integral is less than $
e^{Cp/2} ( \int_{\R^2\times\R^2}f(x,y) \,dx\,dy  )^{p/2} $, which
is equal to $e^{Cp/2}$, hence finite.

\subsection*{Domain $\{\llvert x\rrvert\geq1,\llvert x-y\rrvert\leq1\}$}

The final part of the proof consists in checking that
%
\begin{equation}
\label{lem:M2r} \E \biggl[ \biggl(\int_{\llvert  x\rrvert  \geq1,\llvert  x-y\rrvert  \leq1}f(x,y)M(dx)M(dy)
\biggr)^{p/2} \biggr]<+\infty.
\end{equation}
Because of Lemma~(\ref{lem:f}) item 2, the above relation just boils
down to proving that there exist $\delta>0$ such that
%
\begin{equation}
\label{lem:Mdiaginf} \E \biggl[ \biggl(\mathop{\int_{\llvert  x\rrvert  \geq1}}_{{\llvert  x-y\rrvert  \leq1}}
\frac{\exp
(-\afrac{\llvert  x\rrvert  ^2}{4t} )}{\llvert  x-y\rrvert  ^\delta}M(dx)M(dy) \biggr)^{p/2} \biggr]<+\infty.
\end{equation}

Once again, we first need to estimate the above expectation on stripes
of the type $\{\llvert  x\rrvert  \geq1,2^{-n-1}\leq\llvert  x-y\rrvert  \leq2^{-n}\}$. So we claim
the following.

\begin{lemma}\label{lem:Mdiaginfaux}
Fix $t>0$. For any $\gamma^2<4$ and $p\in\,]1,\frac{4}{\gamma^2}[$, there
exists a constant $C>0$ (only depending on $\E [ M ([0,1]
)^{p} ]$) such that for all $n$:
\begin{eqnarray*}
&& \E \biggl[ \biggl(\mathop{\int_{\llvert  x\rrvert  \geq1}}_{{2^{-n-1}\leq\llvert  x-y\rrvert  \leq2^{-n}}} \exp
\biggl(-\frac{\llvert  x\rrvert  ^2}{4t} \biggr) M(dx)M(dy) \biggr)^{p/2} \biggr]
\leq
C2^{-n(\xi_M(p)-2)}.
\end{eqnarray*}
\end{lemma}

Let us admit for a while the above lemma to finish the proof of~(\ref
{lem:Mdiaginf}). If we choose $\delta$ such that $0<\delta<\frac
{2(\xi
_M(p)-2)}{p}$, we can then use Lemma \ref{lem:Mdiaginfaux} and
sub-additivity to get  to computations~(\ref{om}), similar
\begin{eqnarray*}
&& \E \biggl[ \biggl(\mathop{\int_{\llvert  x\rrvert  \geq1}}_{{\llvert  x-y\rrvert  \leq1}}
\frac
{\exp
(-\afrac{\llvert  x\rrvert  ^2}{4t} )}{\llvert  x-y\rrvert  ^\delta}M(dx)M(dy) \biggr)^{p/2} \biggr] 
\\
&&\qquad \leq\sum_{n=0}^{+\infty}2^{(n+1)\vfrac{\delta p}{2}}C2^{-n(\xi_M(p)-2)},
\end{eqnarray*}
which is a converging series.\vadjust{\goodbreak}

\begin{pf*}{Proof of Lemma \ref{lem:Mdiaginfaux}} We keep on
carrying out the proof in dimension $1$ with the kernel $K(x,y)=\ln
_+\frac{1}{\llvert  x-y\rrvert  }$. It is also plain to check that the expectation is
finite thanks to the exponential term. We will prove the result when
integrating only over the domain $\{x\geq1,2^{-n-1}\leq\llvert  x-y\rrvert  \leq
2^{-n}\}$. It will then be obvious to complete the proof (e.g.,
by using invariance of $M$ in law under reflection). As previously, the
reader may check that the stripe $\{x\geq1,2^{-n-1}\leq\llvert  x-y\rrvert  \leq
2^{-n}\}$ may be covered by the squares $[\frac{k}{2^n},\frac
{k+2}{2^n}]^2$ for $k$ running over the set $K_n=\Z\cap[2^n,+\infty[$.
Therefore, by using subadditivity of the mapping $x\mapsto x^{p/2}$,
we have
\begin{eqnarray*}
&&\E \biggl[ \biggl(\mathop{\int_{x\geq1}}_{{2^{-n-1}\leq\llvert  x-y\rrvert  \leq2^{-n}}} \exp
\biggl(-\frac{\llvert  x\rrvert  ^2}{4t} \biggr) M(dx)M(dy) \biggr)^{p/2} \biggr]
\\
&&\qquad \leq\sum_{k\in K_n} \E
\biggl[ \biggl( \int_{ [\sfrac{k}{2^n},\vfrac
{k+2}{2^n}]^2} \exp \biggl(-\frac{k^2}{ t2^{2n+2}}
\biggr) M(dx)M(dy) \biggr)^{p/2} \biggr]
\\
&&\qquad = \sum_{k\in K_n} \exp \biggl(-\frac{k^2p}{ t2^{2n+2}}
\biggr) \E \biggl[ M \biggl(\biggl[\frac{k}{2^n},\frac{k+2}{2^n}\biggr]
\biggr)^{p} \biggr].
\end{eqnarray*}
By stationarity and scale invariance~(\ref{esi}), we get
\begin{eqnarray*}
&& \sum_{k\in K_n} \exp \biggl(-\frac{k^2p}{ t2^{2n+2}} \biggr)
\E \biggl[ M \biggl(\biggl[\frac{k}{2^n},\frac{k+2}{2^n}\biggr]
\biggr)^{p} \biggr]
\\
&&\qquad    = \sum_{k\in K_n} \exp \biggl(-
\frac{k^2p}{
t2^{2n+2}} \biggr) \E \bigl[ M \bigl(\bigl[0,2^{n-1}\bigr]
\bigr)^{p} \bigr]
\\
&&\qquad    = 2^{-(n-1)\xi_M(p)}\sum_{k\in K_n} \exp
\biggl(-\frac
{k^2p}{ t2^{2n+2}} \biggr)\E \bigl[ M \bigl([0,1] \bigr)^{p}
\bigr]
\leq C2^{-n(\xi_M(p)-1)},
\end{eqnarray*}
where $C$ only depends on $\E [ M ([0,1] )^{p} ]$. The last
line uses the standard trick of convergence of Riemann sums.
\end{pf*}

\begin{rem}\label{rem:mom}
If\vspace*{1pt} $\gamma<2$, it is expected in great generality that $F$ possesses
moments of order $p$ for $p<\frac{4}{\gamma^2}$. We proved that this is
true in the more complicated situation $\sqrt{2}\leq\gamma<2$. If
$0<\gamma<\sqrt{2,}$ we only gave the existence of moments for $p\leq
2$. However, our strategy could be easily adapted to treat the case
$p<\frac{4}{\gamma^2}$. In that case, one has to choose an integer
$n\geq2$ such that $p/n<1$ and apply the Jensen inequality to get an
expression similar to~(\ref{mainmaj}) (replace $2$ by $n$) excepted
that we get an integral over $(\R^2)^n$ instead of $(\R^2)^2$. Then we
can reproduce our strategy up to modifications that are obvious but
notationally awful.
\end{rem}
\end{appendix}

\section*{Acknowledgements}
The authors wish to thank M. Bauer, F.~David, J. Dub\'edat, M.
Gubinelli, J.~F. Le Gall for fruitful discussions and comments. We also
wish to thank the anonymous referee for his/her careful reading.


%
\vspace*{-1pt}

\printaddresses
\end{document}